\documentclass[10pt,a4paper]{amsart}

\usepackage{graphicx}
\usepackage{psfrag}
\include{psfig}

\newcommand{\ncmd}{\newcommand}
\ncmd{\rencmd}{\renewcommand}
\ncmd{\dspst}{\displaystyle }

\ncmd{\tm}{Theorem}

\ncmd{\preu}{\noindent \mbox{\bf Proof  }}
\ncmd{\preud}[1]
{\noindent \mbox{\bf Proof of \ref{#1}   }}

\ncmd{\fin}{\hspace*{\fill} 
\quad\hbox{\hskip 1pt\vrule width 4pt height 6pt
          depth 1.5pt\hskip 1pt} \medskip }

\ncmd{\son}{\mbox{$SO_{0}(1,n-1)\/$ }}
\ncmd{\mink}{\mbox{${\mathbb R}^{1,n-1}\/$ }}
\ncmd{\flt}{\mbox{$\Phi^t\/$ }}

\ncmd{\bas}{\mbox{$\mathcal B\/$ }}
\ncmd{\basd}{\mbox{${\mathcal B}^{\ast}\/$ }}

\newtheorem{prop}{Proposition}[section]
\newtheorem{thm}[prop]{Theorem}
\newtheorem{lem}[prop]{Lemma}
\newtheorem{cor}[prop]{Corollary}
\newtheorem{defin}[prop]{Definition}
\newtheorem{rac}[prop]{Remark}
\newtheorem{ric}[prop]{Example}

\newtheorem{banane}{Figure}

\ncmd{\exe}{\begin{ric} \em }
\ncmd{\eexe}{\em \end{ric}}
\ncmd{\rque}{\begin{rac} \em}
\ncmd{\erque}{\em \end{rac}}

\begin{document}

\title{Three-dimensional Anosov flag manifolds}
\author{Thierry Barbot}
\keywords{Flag manifold, Anosov representation}
\date{May 17, 2005. Work supported 
by CNRS}

\begin{abstract}
Let $\Gamma$ be a surface group of higher genus. 
Let $\rho_0: \Gamma \rightarrow \mbox{PGL}(V)$ be a discrete faithful 
representation with image contained in the natural embedding of 
$\mbox{SL}(2, {\mathbb R})$ in $\mbox{PGL}(3, {\mathbb R})$ as a 
group preserving a point and a disjoint 
projective line in the projective plane. We prove that such a representation is $(G,Y)$-Anosov 
(following the terminology of \cite{labourieanosov}), where $Y$ is 
the frame bundle. More generally, we prove that all the deformations $\rho: \Gamma \rightarrow \mbox{PGL}(3, {\mathbb R})$ studied in \cite{barflag} are $(G,Y)$-Anosov.
As a corollary, we obtain all the main results of \cite{barflag}, 
and extend them to any small deformation of $\rho_0$, not necessarily 
preserving a point or a projective line in the projective space: 
in particular, there is a $\rho(\Gamma)$-invariant solid torus $\Omega$ 
in the flag variety. The quotient space $\rho(\Gamma)\backslash\Omega$ 
is a flag manifold, naturally equipped with two $1$-dimensional 
transversely projective foliations arising from the projections of 
the flag variety on the projective plane and its dual; if $\rho$ is 
strongly irreducible, these foliations are not minimal. More precisely, 
if one of these foliations is minimal, then it is topologically conjugate 
to the strong stable foliation of a double covering of a geodesic flow, 
and $\rho$ preserves a point or a projective line in the projective plane.
All these results hold for any $(G,Y)$-Anosov representation which is not 
quasi-Fuchsian, i.e., does not preserve a strictly convex domain in the projective plane.
\end{abstract}

\maketitle

\section{Introduction}

A flag is a pair $(p,d)$ where $p$ is the point of the projective plane, and $d$ a projective line containing $p$. The group $G = \mbox{PGL}(3,{\mathbb R})$ of projective transformations of 
the projective plane acts naturally on the flag variety $X$, i.e., the space of flags. 

Let $\Gamma$ be the fundamental group of a closed surface $\Sigma$ of 
higher genus. 
In \cite{barflag}, we considered representations 
$\rho: \Gamma \rightarrow G$ near ''horocyclic'' representations, i.e., 
obtained from a faithful discrete representation 
$\Gamma \rightarrow H = \mbox{SL}(2, {\mathbb R})$ 
composed with the natural morphism identifying $H$ 
with the commutator subgroup of the stabilizer in $G$ of 
a point $p_0$ and a projective line $d_0$ in the projective plane, 
with $p_0 \notin d_0$. Actually, in \cite{barflag}, we only considered 
some deformations of horocyclic representations for which $p_0$ is still 
a global fixed point: we called these representations 
\emph{hyperbolic representations.\/} We proved that for such a 
representation, there is a closed $\rho(\Gamma)$-invariant 
simple closed curve $\Lambda$ in $X$, and a open $\rho(\Gamma)$-invariant 
domain $\Omega$, both depending on $\rho$, such that:

-- $\Lambda$ is the image of a $1$ to $1$ continuous $\Gamma$-equivariant map from the projective line ${\mathbb R}P^1$ into $X$, where the action of $\Gamma$ on ${\mathbb R}P^1$ is the usual projective action (which is unique up to topological conjugacy).

-- The action of $\rho(\Gamma)$ on $\Omega$ is free and properly 
discontinuous. The quotient space of this action, called $M$, is 
a flag manifold (cf. \S~\ref{deflagstru}). The first tautological 
foliation (see \S~\ref{tautogoldman}) is topologically conjugate 
to the strong stable foliation of a double covering of the geodesic 
flow on $\Sigma$. On the other hand, the second tautological 
foliation is not minimal, except when $\rho(\Gamma)$ also preserves 
a projective line.

In the present paper, we extend all these results omitting 
the assumption that $\rho(\Gamma)$ admits a global fixed point.
The key observation is that horocyclic representations 
are $(G,Y)$-Anosov in the terminology of \cite{labourieanosov}, 
where $Y$ is the frame variety, i.e., the space of non-collinear 
points in the projective plane (see \S~\ref{deframe}).

Typical $(G,Y)$-Anosov representations are \emph{hyperconvex,\/} 
i.e., those preserving a strictly convex domain of the 
projective plane. S. Choi and W. Goldman proved that every representation
$\rho: \Gamma \rightarrow G$ which is quasi-Fuchsian, i.e., which can 
be continuously deformed to a representation taking values 
in $\mbox{SO}_0(1,2) \subset \mbox{PGL}(3,{\mathbb R})$, preserves such a 
strictly convex curve (\cite{choigoldman}). It follows easily that they 
are $(G,Y)$-Anosov, hence, hyperconvex. F. Labourie has extended this result 
to the higher-dimensional case (\cite{labourieanosov}). In \cite{guichard}, 
O. Guichard proved that conversely, any hyperconvex representation is 
quasi-Fuchsian.

Here, we consider general $(G,Y)$-Anosov representations, 
not necessarily hyperconvex. Each preserves a 
limit curve $\Lambda$, which is a H\"{o}lder continuous simple closed curve 
(see definition~\ref{deflambda}). In \S~\ref{general}, we establish some general results
on this limit curve. In addition, we prove that such a representation 
always preserves an open domain $\Omega$ on which it acts 
freely and properly discontinuously (Theorem~\ref{actiopropre}). We insist 
here on the low regularity of $\Lambda$ and of the boundary $\partial\Omega$: 
when $\rho$ is irreducible and not hyperconvex, then they are not Lipschitz regular 
(see corollary~\ref{pasregulier}). Let's mention here that a similar statement 
is true in the hyperconvex case: the limit curve is in general $C^1$ with H\"older 
derivatives, but if the derivatives are Lipschitz, then $\rho(\Gamma)$ 
preserves a conic, i.e., is conjugate in $\mbox{PGL}(3, {\mathbb R})$ to 
a Fuchsian subgroup.

Our first interest is in non-hyperconvex $(G,Y)$-Anosov 
representations: we prove that hyperbolic representations 
are $(G,Y)$-Anosov (Theorem~\ref{flaganosov}, Lemma~\ref{chouia}). 
We actually suspect that non-hyperconvex $(G,Y)$-Anosov representations 
form a connected space (see Question $2$ in \S~\ref{conclu}). Observe 
that they all belong to the same connected component of the space of 
representations of $\Gamma$ into $G$: the component of the trivial 
representation (corollary~\ref{comptriviale}), whereas hyperconvex 
representations are those in the Hitchin component (see \S~\ref{espacedef}, 
Remark~\ref{qfs}).

For any $(G,Y)$-Anosov representation $\rho$, the quotient 
$M = \rho(\Gamma)\backslash\Omega$ is a natural flag manifold, 
which we call an \emph{Anosov flag manifold.\/} It is therefore naturally 
equipped with two transversely projective $1$-dimensional foliations: 
the tautological foliations (see \S~\ref{tautogoldman}). When $\rho$ is 
hyperconvex, $M$ admits three connected components, with well-understood 
geometrical features (see remark~\ref{fauxlorentz}), and 
tautological foliations with well-identified dynamical properties (see 
Remark~\ref{tautoFuchsian}): they are either foliations by circles
or doubly covered by the geodesic flow on $\Sigma$ for any hyperbolic metric.

The situation when $\rho$ is not hyperconvex is completely different 
(cf. \S~\ref{pashyper}): the tautological foliations in this case 
are never foliations by circles or finitely covered by Anosov flows. 
The correct picture is hard to capture: in some cases, they are topologically 
conjugate to the horocyclic flow of some Anosov flow 
(see Proposition~\ref{conjuguons}), but this is not true in general. 
For example, when $\rho$ is strongly irreducible, the tautological 
foliations are \emph{not\/} minimal, in contrast to the horocyclic flow. 
The dynamical properties of these foliations are quite interesting.
We suspect that these tautological foliations never admit periodic orbits 
(see Question $5$ in \S~\ref{conclu}). If our suspicion is confirmed, 
it would provide examples of flows with unusual
behavior. For example, recall the Seifert Conjecture, asserting that any flow
on the three-dimensional sphere admits a periodic orbit. The first smooth
counterexamples to this conjecture were found by K. Kuperberg (\cite{kuper1}).
Observe moreover that the tautological foliations
considered here can be volume-preserving: for example, this is the case for 
second tautological foliations associated to hyperbolic 
radial representations for which the morphism $u: \Gamma \rightarrow {\mathbb R}$
is trivial (see \S~\ref{flagoldman}). As far as we know, the only known 
examples of volume preserving flows on the $3$-sphere without periodic orbits
have regularity at most $C^2$ (\cite{kuper2}, \cite{ginzburg}). As a matter
of fact, volume preserving flows on $3$-manifolds which are not minimal and
without periodic orbits are quite uncommon; hence, it seems to us quite 
interesting to answer our Question $5$.

\section*{Acknowledgements}
I would like to thank F. Labourie, who made me aware of the $(G,Y)$-Anosov 
character of horocyclic representations. G. Kuperberg has corrected 
some imprecisions in a previous version of this paper regarding counterexamples 
to the Seifert conjecture.

\section{Definitions}

\subsection{The Flag variety}
\label{deflag}

-- Let $V$ be the vector space ${\mathbb R}^3$, and $(e_1, e_2, e_3)$ 
its canonical basis. Let $V^\ast$ be the dual vector space, with the dual canonical basis
$(e^\ast_1, e^{\ast}_2, e^{\ast}_{3})$. We denote by $\langle v \mid
v^\ast \rangle$ the evaluation of an element $v^\ast$ of $V^{\ast}$ on 
an element $v$ of $V$.

-- Let $N(v)$, $N(v^\ast)$ denote the norms on $V$, $V^\ast$, respectively, 
for the Euclidean metrics on $V$, $V^\ast$ in which the canonical 
basis is orthonormal.

-- $P(V)$ and $P(V^{\ast})$ are the associated projective spaces. Elements of
$P(V)$ are denoted $[v]$.

-- The \emph{flag variety\/} $X$ is the closed subset of $P(V) \times P(V^{\ast})$ formed by pairs $([v], [v^{\ast}])$ satisfying $\langle v \mid v^\ast \rangle = 0$.  

-- $G$ is the group $\mbox{SL}(V) \approx \mbox{PGL}(V)$.
The group $G$ acts naturally on $V$ and admits a dual (left) action on 
$V^{\ast}$ uniquely defined by requiring 
$\langle u \mid g.v^\ast \rangle = \langle g^{-1}u \mid v^\ast \rangle$ 
for any  $u$, $v^\ast$. For any $g$ in $G = \mbox{SL}(V)$, we 
denote by $g^\ast$ the corresponding element of $\mbox{SL}(V^\ast)$. 
If $\mbox{SL}(V)$ and $\mbox{SL}(V^\ast)$ are identified with 
$\mbox{SL}(3, {\mathbb R})$ via the canonical basis, 
$g^\ast$ is the inverse of the transpose of $g$.

-- The diagonal action restricts as a natural action of $G$ on ${X}$.

\rque
\label{precision}
Every element of $P(V^\ast)$ defines a projective line in $P(V)$
(the projection to $P(V)$ of its kernel). Dually, every element of
$P(V^\ast)$ corresponds to a projective line in $P(V)$. Hence, we 
can consider $P(V)$ to be the space of projective lines of $P(V^{\ast})$ 
and $P(V^{\ast})$
to be the space of projective lines of $P(V)$. 

The image we have in mind is to view an element of $X$ as a flag, 
i.e., a point in the projective plane and a projective line containing 
this point. 

In order to formalize this point of view, we introduce 
the following notation: if $K$ is a subspace of $V$ (resp. of $V^\ast$), we 
denote by $K^\perp$ its orthogonal, i.e., the subset of $V^\ast$ 
(resp. of $V$) vanishing on $K$. Hence, 
for any $[v]$ in $P(V)$, $[v^\perp]$ is the corresponding 
projective line in $P(V^\ast)$.

If $[u]$, $[v]$ are distinct elements of $P(V)$, we denote by 
$[(uv)^\ast]$ the element $[K^\perp]$ of $P(V^\ast)$, where 
$K$ is the $2$-dimensional space spanned by $u$ and $v$. We employ 
similar notation when $[u]$, $[v]$ belong to $P(V^\ast)$.
\erque

\rque
\label{identifions}
Let $P^+$ (respectively $P^{-})$ be the subgroup of $G$ containing the upper 
(respectively lower) triangular matrices 
(under the identification of $G$ with $\mbox{SL}(3, {\mathbb R})$). 
Let $X^{\pm}$ be the quotient spaces $G/P^{\pm}$:  the map 
$g \mapsto ([ge_1], [g^{\ast}e^\ast_3])$ induces the 
identification $X^+ = G/P^+ \approx X$, and the map 
$g \mapsto ([ge_3], [g^{\ast}e^\ast_1])$ induces the 
identification $X^- = G/P^- \approx X$. These 
identifications are $G$-equivariant.
\erque

\subsection{Flag manifolds}
\label{deflagstru}
A flag structure is a $(G,X)$-structure, for $(G,X)$ as above. We briefly present this notion here. For a more complete description of $(G,X)$-structures, see \cite{thurston} or \cite{goldman}.

A flag structure on a manifold $M$  is an atlas on $M$ with charts 
taking values in $X$ and coordinates changes expressed in the 
charts by restrictions of elements of $G$. A typical example of a 
flag structure is the quotient of an open domain $\Omega$ of $X$ by a 
discrete subgroup of $G$ acting freely and properly discontinuously on 
$\Omega$; in particular, the flag variety $X$, itself, is an example. 
Less trivial examples are given in \S~\ref{exempleflagmanifold}, 
\S~\ref{anosovmanifold}.

A flag map between flag manifolds is a map which can be locally 
expressed in the flag charts by restrictions of elements of $G$. 
A flag map is always a local homeomorphism. A flag map which is bijective 
is called a flag isomorphism. A flag manifold is an isomorphism class of 
flag structures on the manifold.

Let $p: \widetilde{M} \rightarrow M$ be the universal covering and 
$\widetilde{\Gamma}$ the fundamental group of $M$, viewed as 
the group of covering automorphisms of $p$. For any flag structure 
on $M$, there is a map ${\mathcal D}: \widetilde{M} \rightarrow X$, called 
the \emph{developping map,\/} and a representation 
$\rho: \widetilde{\Gamma} \rightarrow G$, called 
\emph{holonomy representation,\/} such that:

-- the maps $p$ and $\mathcal D$ are flag maps,

-- $\mathcal D$ is $\widetilde{\Gamma}$-equivariant: 

\[ \forall \gamma \in \widetilde{\Gamma} \;\; {\mathcal D} \circ \gamma = \rho(\gamma) \circ \mathcal D \]

\subsection{The frame variety}
\label{deframe}

Let $Y$ denote the \emph{frame variety}, i.e., the space of triples $([u], [v], [w])$ of noncollinear elements of $P(V)$. This $6$-dimensional space is homogeneous under the diagonal action of $G$: it can be identified with the right quotient $G / Z$, where $Z$ is the group of diagonal matrices. 

The frame variety admits two natural projections $\pi_\pm$ on $X$, 
defined by (see remark~\ref{precision} for the notation):

\[ \pi_+([u], [v], [w]) = ([u], [(uv)^\ast]) \;\; ; \; \pi_-([u], [v], [w])= ([w], [(wv)^\ast]) \]

These projections are both $G$-equivariant. 
Together, they define a map $\Pi: Y \rightarrow {X} \times {X}$. 
We denote by $\mathcal Y$ the image of $\Pi$. 

\begin{lem}
\label{defY}
The elements of $\mathcal Y$ are the 
pairs $( ([u], [u^\ast]), ([v],[v^\ast]))$ satisfying:

-- $\langle u \mid u^\ast \rangle = 0$,

-- $\langle v \mid v^\ast \rangle = 0$,

-- $\langle u \mid v^\ast \rangle \neq 0$,

-- $\langle v \mid u^\ast \rangle \neq 0$.\fin
\end{lem}

In other words, $[u]$ (respectively $[v]$) must belong to 
the projective line $[u^\ast]$ (resp. $[v^\ast]$) since we 
are considering elements of $X$, but it cannot belong to 
the projective line $[v^\ast]$ (respectively $[u^\ast]$).

The fibers of the projections $\pi_\pm$ are the leaves 
of $G$-invariant foliations ${\mathcal G}^{\pm}$ of $X$. 
Let $E^\pm$ be the tangent bundles to these foliations. We obtain a 
$G$-invariant decomposition $TY = E^+ \oplus E^-$ of the tangent bundle of $Y$.


\rque
The inclusions $Z \subset P^\pm$ defines canonical maps 
$G/Z \rightarrow X^{\pm}$, which,  via the identifications 
$X^\pm \approx X$ presented in remark~\ref{identifions} and $G/Z \approx Y$, 
are the maps $\pi_\pm$. 
\erque

\subsection{The geodesic flow as an Anosov flow}
\label{geodanosov}
Let $\Sigma$ be a closed surface with negative Euler characteristic, 
and $\Gamma$ the fundamental group of $\Sigma$. Select a $1$ to $1$ morphism 
${\i}: \Gamma \rightarrow H$ with discrete image, where $H$ denotes the group 
$\mbox{SL}(2, {\mathbb R})$. This induces a Fuchsian representation 
${\bar{\i}}: \Gamma \rightarrow \overline{H}$, where 
$\overline{H}$ denotes the group $\mbox{PSL}(2, {\mathbb R})$ 
(observe that conversely, any Fuchsian representation lifts to 
a representation into $H$, since the associated Euler class is even). 

Consider the flow on $H$ induced by the right action of the $1$-parameter 
group $A = \{ a^t \}$, where:

\[ a^{t} = \left(\begin{array}{cc}
           e^t & 0 \\
	   0 & e^{-t}
	   \end{array}\right) \]

This flow induces another flow on $M = {\i}(\Gamma) \backslash H$, 
denoted by $\Phi^t$. Its projection on 
$\overline{M} = {\bar{\i}}(\Gamma) \backslash \overline{H} \approx T^1\Sigma$ 
is a flow $\overline{\Phi}^t$.

\rque
\label{phigeod}
This flow divided by $2$, i.e., the flow $p \mapsto \Phi^{t/2}(p)$, is the familiar geodesic flow
 associated to the Riemannian surface ${\bar{\i}}(\Gamma) \backslash {\mathbb H}^2 \approx \Sigma$, where ${\mathbb H}^2$ is the hyperbolic plane. 
\erque

Consider the following $1$-parameter subgroups of $H$:

\[    h^{s}_{+} = \left(\begin{array}{cc}
           1 & s \\
	   0 & 1
	   \end{array}\right) \; ; \; 
	   h^{s}_{-} = \left(\begin{array}{cc}
           1 & 0 \\
	   s & 1
	   \end{array}\right) \]

We have the identities:
\begin{eqnarray} 
h^s_+ a^t= a^t h^{\exp(-2t)s}_+\; ; \; h^s_- a^t= a^t h^{\exp(2t)s}_- 
\end{eqnarray}

Denote also by $h^s_{\pm}$ the flow on $M$ induced by the right actions 
of these $1$-parameter subgroups. The flow  $h^{s}_{+}$ (resp. $h^s_-$) 
is called the \emph{stable horocyclic flow\/} 
(resp. \emph{unstable horocyclic flow\/}). 
The flow $\Phi^t$ permutes the orbits of $h^s_{\pm}$. 
Moreover, the orbits of $h^s_+$ are exponentially contracted by $\Phi^t$, 
whereas the orbits of $h^s_-$ are exponentially expanded by $\Phi^t$. 
This feature establishes precisely that $\Phi^t$ and $\overline{\Phi}^t$ are \emph{Anosov.\/}

\begin{defin}
\label{defanosov}
An Anosov flow on a closed manifold $M$ equipped with a 
Riemannian metric $\Vert$ is a non-singular flow $\Phi^t$ 
such that the differential of $\Phi^t$ preserves a decomposition 
$TM = \Delta \oplus E^{ss} \oplus E^{uu}$ of the tangent bundle, 
 satisfying the following properties, for some positive constants $a$, $b$:

\begin{itemize}
\item The line bundle $\Delta$ is tangent to the flow,

\item for any vector $v$ in $E^{ss}$ over a point $x$ of $M$, and for any positive $t$:

\[ \Vert D_{x}\Phi^{t}(v) \Vert \leq b e^{-at} \Vert v \Vert \]

\item for any vector $v$ in $E^{uu}$ over a point $x$ of $M$, and for any negative $t$:

\[ \Vert D_{x}\Phi^{t}(v) \Vert \leq b e^{at} \Vert v \Vert \]

\end{itemize}
\end{defin}

\subsubsection{Stable and unstable leaf spaces}

Denote by ${\mathcal A}_{\pm}$ the group generated by $a^t$ and $h^s_{\pm}$.
It is isomorphic to the group of volume preserving affine transformations 
of the plane. The orbits on $M$ or $\overline{M}$ of ${\mathcal A}_+$ 
are called \emph{stable leaves;\/} the orbits of ${\mathcal A}_-$ 
are called \emph{unstable leaves.\/} We denote by ${\mathcal S}_\pm$, 
$\overline{\mathcal S}_{\pm}$ the quotient spaces $H / {A_\pm}$,
$\overline{H} / {\mathcal A}_\pm$. The latter, $\overline{\mathcal S}_\pm$, 
are both homeomorphic to the projective line ${\mathbb R}P^1$, and 
${\mathcal S}_\pm$ are double coverings of $\overline{\mathcal S}_{\pm}$. 
Moreover, these identifications are $H$-equivariant, where the $H$-action 
on $\overline{\mathcal S}_\pm$ is the action induced by left translation, 
and the $H$-action on ${\mathbb R}P^1$ is the usual projective action.

\subsubsection{The bifoliated orbit space}
\label{bifoliated}
The map $hA \mapsto (h{\mathcal A}_+, h{\mathcal A}_-)$ embeds the orbit 
space $\overline{Q} = \overline{H} / a^t$ into the torus 
$\overline{\mathcal S}_{+} \times \overline{\mathcal S}_-$. 
More precisely, for fixed $H$-equivariant identifications
$\overline{\mathcal S}_\pm \approx {\mathbb R}P^1$, the image of 
this embedding is the complement in ${\mathbb R}P^1 \times {\mathbb R}P^1$ 
of the diagonal $\Delta$.
In other words, every ${\mathcal A}_+$-orbit $x$ intersects every
${\mathcal A}_-$-orbit, \emph{except one,\/} which we call $\alpha(x)$. 
Note that we have defined a continuous $\Gamma$-equivariant 
map $\alpha: \overline{\mathcal S}_+ \rightarrow \overline{\mathcal S}_-$. 

We denote by $\overline{\mathcal Q}$ the image of $\overline{Q}$ in 
$\overline{\mathcal S}_{+} \times \overline{\mathcal S}_-$; 
this is the complement of the graph of 
$\alpha: \overline{\mathcal S}_+ \rightarrow \overline{\mathcal S}_-$.



\subsection{Anosov representations}
\label{zalabourie}

Let $\rho: \Gamma \rightarrow G$ be any representation, and 
let $\pi_\rho: E_\rho \rightarrow M$ be the associated 
flat $(G,Y)$-bundle: $E_\rho$ is the quotient of 
$H \times Y$ by the relation identifying each $(h, y)$ with 
$({\i}(\gamma) h , \rho(\gamma).y)$, for every $\gamma$ in 
$\Gamma$. The projection $(h,y) \mapsto h$ induces a map 
$\pi_\rho$ from $E_\rho$ onto $M = {\i}(\Gamma)\backslash H$, 
which is a $G$-bundle, with fiber $Y$.

The (trivial) foliation of $H \times Y$ having as leaves the fibers 
of $(h, y) \mapsto y$ induces a foliation on $M$, which we 
denote by ${\mathcal F}_\rho$, and call the \emph{horizontal foliation.\/} 
The leaves of ${\mathcal F}_\rho$ are transverse to the fibers of $\pi_\rho$.

The flow $\Phi^t$ lifts uniquely to a horizontal flow ${\Phi}_\rho^t$ in 
$E_\rho$, i.e., tangent to the horizontal foliation: 
just take the flow induced in the quotient by the flow 
on $H \times Y$ defined by $(h, y) \mapsto (ha^t, y)$.

We have defined foliations ${\mathcal G}^\pm$ on $Y$. They provide 
two $3$-dimensional foliations on $H \times Y$ and induce 
on $E_\rho$ two $3$-dimensional foliations ${\mathcal F}^\pm$ 
which are preserved by $\Phi^t_\rho$ and tangent to the fibers of $\pi_\rho$.

We will mainly consider the tangent bundles $E_\rho^\pm$ of 
these foliations, which are canonically induced by the bundles 
$E^{\pm}$ over $Y$.

\begin{defin}[\cite{labourieanosov}]
\label{defrepanosov}
A $(G,Y)$-Anosov structure over $(M, \Phi^{t})$ is the data 
of a representation $\rho: \Gamma \rightarrow G$ and 
a continuous section $s$ of $\pi_\rho$ with the following properties:

-- the image $S$ of $s$ is $\Phi_\rho^t$-invariant,

--  for any norm $\Vert$ on $E_\rho$, there are positive constants $a$ and $b$ 
such that for any $p$ in $S$, any positive $t$, any vector $v^+$ 
in $E^{+}_\rho$ over $p$, and any vector $v^-$ in $E^-_\rho$ over $p$, we have:
\[ \Vert D_{p}\Phi_\rho^{t}(v^+) \Vert \leq b e^{-at} \Vert v^+ \Vert \]
\[ \Vert D_{p}\Phi_\rho^{-t}(v^-) \Vert \leq b e^{-at} \Vert v^- \Vert \]

If these conditions are fulfilled, then the representation $\rho$ is said to be 
$(G,Y)$-Anosov, or an \emph{Anosov flag representation.\/}
\end{defin}

\rque
Here,  $Y$ is the frame variety, and $G$ the group $\mbox{SL}(V)$, but, 
of course, the definition~\ref{defrepanosov} extends to any other pair $(G, Y)$, and 
any Anosov flow can play the role played here by the ''geodesic flow'' 
$\Phi^t$.
\erque

The section $s$ appearing in this definition is not assumed to be 
differentiable, even nor Lipschitz. The maximal regularity which can 
be required in general is H\"{o}lder
continuity. 

The main interest of Anosov representations is their stability: it follows from the structural stability of Anosov flows (more precisely, of hyperbolic closed sets) that the set of $(G,Y)$-Anosov representations is an open domain
in the space of representations of $\Gamma$ into $G$ equipped with its natural 
topology (see proposition $2.1$ of \cite{labourieanosov}). Observe that 
the perturbed $(G,Y)$-Anosov structure covers the \emph{same\/} Anosov flow $(M, \Phi^t)$.

Another important feature of Anosov representations is that they provide nice $\Gamma$-invariant geometric objects, obtained as follows:

The section $s$ of a $(G,Y)$-Anosov structure lifts to a continuous map $f: H \rightarrow Y$ such that:

-- $f$ is $\Gamma$-equivariant: $ f \circ {\i}(\gamma) = \rho(\gamma) \circ f$,

-- $f$ is invariant by the lifted flow: $f(ha^t) = f(h)$.
 
Therefore, $f$ induces a continuous $\Gamma$-equivariant map from the orbit space $Q = H / A$ into $Y$.

\begin{lem}
\label{pif}
The restriction of $\pi_+ \circ f$ (resp. $\pi_- \circ f$) to any ${\mathcal A}_+$-orbit (resp. ${\mathcal A}_-$-orbit) is constant.
\end{lem}

\preu
Let $p$, $p'$ be two elements of $M$ belonging to the same stable leaf, 
i.e., on the same (right) ${\mathcal A}_+$-orbit: 
for positive times $t$, the iterates $\Phi^{t}(p)$ and $\Phi^{t}(p')$ remain 
at a bounded distance apart inside the stable leaf.
Then the exponential dilatation along the leaves of ${\mathcal F}^{-}$ implies that $s(p)$ and $s(p')$ must belong to the same leaf of ${\mathcal F}^{+}$. The lemma follows.\fin

It follows from the lemma~\ref{pif} that $\pi_+ \circ f$ and $\pi_- \circ f$ 
induce maps
$f_+: {\mathcal S}_+  \rightarrow {X}$ and 
$f_-: {\mathcal S}_- \rightarrow {X}$. Of course, 
$f_\pm$ are both $\Gamma$-equivariant: their images (which, as we will see, 
are the same) are $\Gamma$-invariant (\emph{a priori\/} only immersed) 
topological circles in the flag variety.

\subsubsection{Splitting the definition}

Consider the flat $X$-bundle $\pi_{\rho}^X: E_{\rho}(X) \rightarrow M$ 
associated to $\rho$. The total space $E_{\rho}(X)$ is the quotient of 
$H \times X$ by the $\Gamma$-action defined by 
$(h,x) \rightarrow (\gamma h, \rho(\gamma)x)$. 
The map $\Pi$ defines a fibered embedding 
$\Pi_{\rho}: E_{\rho} \rightarrow E_{\rho}(X) \times E_{\rho}(X)$. Let 
${\mathcal Y}_\rho$ be the image of $\Pi_\rho$. The flow $\Phi^t$ also 
lifts in a unique way to a horizontal flow $\Phi_X^t$ on $E_{\rho}(X)$: 
$\Phi^t_X(h, x) = (ha^t, x)$. Obviously, Definition~\ref{defrepanosov}
is equivalent to:

\begin{defin}
\label{def2rep}
A $(G,Y)$-Anosov structure over $(M, \Phi^{t})$ is the data of a representation $\rho: \Gamma \rightarrow G$ and two continuous section $s^\pm$ of $\pi_\rho^X$ satisfying the following properties:

-- the sections $s^\pm$ are preserved by the flows: $\Phi^t_X(s^\pm(p)) = s^\pm(\Phi^tp)$,

--  the image of $s^+$ (resp. $s^-$) is a (exponentially) repellor (resp. attractor) for $\Phi^t_X$,

-- for every $p$ in $M$, the pair $(s_+(p), s_-(p))$ belongs to ${\mathcal Y}_{\rho}$.

\end{defin}

\subsubsection{Quasi-Fuchsian representations}
\label{qfs}
A very nice family of $(G,Y)$-Anosov representations is the 
family of \emph{quasi-Fuchsian\/} representations 
(in the terminology of \cite{labourieanosov}), i.e, 
the representations 
$\rho: \Gamma \rightarrow G  \approx \mbox{SL}(3, {\mathbb R})$ 
which are in the \emph{Hitchin component,\/} i.e., which can be 
deformed to \emph{Fuchsian representations,\/} in other words, 
to a representation $\rho_0: \Gamma \rightarrow \mbox{SO}_0(1,2)$. 
Indeed, F. Labourie has proven that quasi-Fuchsian representations 
are $(G,Y)$-Anosov. 

Let's be a bit more precise: S. Choi and B. Goldman proved in 
\cite{choigoldman} that any representation in the Hitchin component 
induces an action on $P(V)$ preserving a strictly convex domain $C$. 
The set of flags $([u], [v^\ast])$ where $[u]$ is a point in $\partial C$ 
and $[v^\ast]$ a support projective line of $C$ is a $\rho(\Gamma)$-invariant 
curve in ${X}$; in fact, this topological circle is equal 
to the images $f_+({\mathcal S}_{+})$ and $f_-({\mathcal S}_{-})$. 
In order to provide a simple hint, we just claim that the fact that 
this kind of representation is flag Anosov can be inferred from the fact 
that the geodesic flow associated to the Finsler Hilbert metric on 
$\rho(\Gamma)\backslash C$ is Anosov.

Following F. Labourie, this kind of curve is said to be \emph{hyperconvex,\/} and
a $(G,Y)$-Anosov representation preserving a hyperconvex curve is said 
hyperconvex. Thus, quasi-Fuchsian representations and hyperconvex representations
coincide. O. Guichard extended recently this statement for all dimensions 
(\cite{guichard}). 

\subsubsection{Canonical Anosov flag representations}

Here, we consider another family of $(G,Y)$-Anosov representations, 
which are not quasi-Fuchsian. They are obtained from an embedding $\rho_0$ of 
$H = \mbox{SL}(2, {\mathbb R})$ into a subgroup of 
$G = \mbox{SL}(V)$ that admits a global fixed point in $P(V) \times P(V^{\ast}) \setminus X$. 

For the fixed point, we select here the pair $(e_2, e^{\ast}_{2})$. The embedding $\rho_0$ is the representation sending the matrix:

\[ \left(\begin{array}{cc}
         a & b \\
	 c & d
	 \end{array}\right) \]
	 
to the element 

\[ \left(\begin{array}{ccc}
        a & 0 & b \\
	0 & 1 & 0 \\
	c & 0 & d 
	\end{array}\right) \]
	
of $G \approx \mbox{SL}(3, {\mathbb R})$ (in the identification arising from the canonical basis).

Recall the notation introduced in \S~\ref{deflag} and \ref{deframe}. Observe 
that 
$\rho_{0}(A)$ is contained in $Z$, and that every $\rho_{0}(h_{\pm}^{s})$ 
belongs to $P^{\pm}$. The adjoint action of $\rho_{0}(a^t)$ on the 
Lie algebra $sl(V)$ of $G$ is diagonalizable, and our choice of  $\rho_0$ 
ensures that the subspace spanned by eigenvectors 
with positive eigenvalues (resp. negative) for $ad(a^{t})$ is the Lie algebra ${\mathcal P}^{+}$ (resp. ${\mathcal P}^{-}$) of $P^+$ (resp. $P^-$). Their intersection is the Lie algeba of $Z$.

\begin{thm}[Proposition $3.1$ of \cite{labourieanosov}]
\label{cestanosov}
Let $\bar{\Gamma}$ be any cocompact subgroup of $H$. The restriction of $\rho_{0}$ to $\bar{\Gamma}$ is $(G,Y)$-Anosov.
\end{thm}

\preu
See Appendix A. \fin

\rque
\label{uniform}
Of course, we will applyTheorem~\ref{cestanosov} to 
$\bar{\Gamma} = {\i}(\Gamma)$. Actually, from the beginning,
we could have selected as $\Gamma$ any discrete cocompact 
subgroup of $H$, possibly with torsion, but this level of 
generality requires a little more caution in the formulation of statements, 
which we considered unnecessary and slightly uncomfortable. 
The reader should have no difficulty to extend the results of this paper 
to this more general context.
\erque

\begin{defin}
The $(G,Y)$-Anosov representation $\rho_0 \circ \i$ is a \emph{canonical Anosov flag representation.}
\end{defin}

\rque
Henceforth, except in Proposition~\ref{rigid}, we will drop 
the symbol $\i$, considering $\Gamma$ directly as a discrete subgroup of $H$.
\erque

\subsubsection{Invariant curves for canonical Anosov representations}
Consider a cano\-nical Anosov flag representation $\rho_0: \Gamma \rightarrow G$.
According to lemma~\ref{pif}, there are two $\Gamma$-equivariant maps $f_\pm:
{\mathcal S}_\pm \rightarrow X$, coming from a $\Gamma$-equivariant map
$f: H \rightarrow Y$. Here, the map $f$  is defined by: $f(h) =  ([\rho_{0}(h)e_1],
[e_2], [\rho_{0}(h)e_3])$ (see Appendix A). 
Hence, we have: $f_+(h{\mathcal A}_{+})  = ([\rho_{0}(h)e_1],
[\rho_{0}(h)^{\ast}e_3^\ast])$, and $f_{-}(h{\mathcal A}_-) = ([\rho_{0}(h)e_3],
[\rho_{0}(h)^{\ast}e_1^{\ast}])$. These two maps have the same image
$\Lambda_0$, which is the set of pairs $([u], [u^{\ast}])$ where $[u]$
belongs to the projective line $[(e^\ast_2)^\perp]$ and $[u^\ast]$ is a 
projective line containing $[e_2]$ (in a more symmetric formulation, $[u^\ast]$ belongs to the projective line 
$[e_2^\perp]$ in $P(V^\ast)$).

We denote by $L_0$ the projective line
$[(e^\ast_{2})^\perp]$ of $P(V)$, and by 
$L^\ast_0$ the projective line $[e^\perp_2]$ of $P(V^\ast)$. We can then 
reformulate the statement above: \emph{a point $([u], [u^\ast])$ in $X$ 
belongs to the curve $\Lambda_0$ if and only if 
$[u]$ belongs to $L_{0}$ and $[u^\ast]$ belongs to $L_{0}^{\ast}$.\/}
Observe that for every  $[u]$ in $L_0$ there is one and only one element 
$[u^\ast]$ of $L_{0}^\ast$ containing $[u]$: the projective line containing 
$[u]$ and $[e_2]$. 

Observe also that \emph{$\Lambda_0$ is the closure of the set of 
attractive fixed points of elements of
$\rho_0({\Gamma})$ in $X$.\/}

Finally, \emph{the image of $f$ is the space of triples 
$([u], [e_2], [w])$ of elements of $P(V)$ where $[u]$ and $[w]$ are distinct elements of $L_0$.\/}

\subsubsection{Canonical Anosov flag manifolds.}
\label{exempleflagmanifold}
Consider the orbits of $\rho_0(H)$
in $X$: there is one $1$-dimensional orbit, the curve $\Lambda_0$. 
There are two $2$-dimensional orbits: 

-- the orbit $A_0$ containing the points $([u], [u^\ast])$ where $[u]$ belongs to $L_0$ and
$[u^\ast]$ does not belong to $L_0^\ast$,

-- the orbit $A^\ast_0$ containing the points $([u], [u^\ast])$ where $[u]$ does not belong to $L_0$ and $[u^\ast]$ belongs to $L_0^\ast$,

There is one open orbit: $\Omega_0 = \{ ([u], [u^\ast]) \in {X} \;/ \; [u] \notin L_0, \;[u^\ast] \notin L^\ast_0 \}$.

The closures $T_0$, $T^\ast_0$ of $A_0$, $A^\ast_0$ are Klein bottles, the intersection $T_0 \cap T^\ast_0$ is $\Lambda_0$, and $\Omega_0$ is the complement in $X$ of $T_0 \cup T^\ast_0$.

The action of $\rho_{0}(H)$ on $\Omega_0$ is simply transitive. It provides 
an identification
$\Omega_0 \approx H$. Therefore, the manifold $M = \Gamma \backslash H$ is 
homeomorphic to the quotient of 
$\Omega_0$ by $\rho_0(\Gamma)$. It provides a natural flag structure on $M$.

\begin{defin}
\label{defcan}
The quotient $\rho_0(\Gamma)\backslash\Omega_0$ is a canonical Anosov flag manifold.
\end{defin}

\subsection{Existence of deformations}
\label{espacedef}

Consider the space $\mbox{Rep}(\Gamma, G)$ of representations of the surface group into $G$, modulo inner automorphisms of $G$ on the target.

\begin{thm}[Hitchin]
$\mbox{Rep}(\Gamma, G)$ has three connected components.
\end{thm}
Let's briefly discuss each of these $3$ components:

-- \emph{The Hitchin component:\/} this is the component containing the Fuchsian representations. Elements of this component are represented by quasi-Fuchsian representations (see \S~\ref{qfs}).

--  \emph{The trivial component:\/} this is the component containing the trivial representation.

-- Representations in the third component are characterized by the fact that they do not lift to
representations from $\Gamma$ into the double covering  $P^{+}\mbox{GL}(V)$
of $G$.

Canonical Anosov flag representations are not quasi-Fuchsian, and they clearly lift to representations in $\mbox{GL}(V)$: hence, they belong to the trivial component. It follows immediately that they can be deformed to representations in $G$ which are not canonical! But there are much more elementary ways to prove this statement: any 
canonical representation $\rho_0$ can be deformed to a strongly irreducible representation, i.e., with image containing no finite index subgroup stabilizing a point or a projective line in $P(V)$
(see e.g. Proposition $3.11$ of \cite{barflag}).

\section{General properties of $(G,Y)$-Anosov representations}
\label{general}
Let $\rho: \Gamma \rightarrow G$ be any $(G,Y)$-Anosov representation. 
According to Lemma~\ref{pif}, there is a $G$-equivariant map 
$f: H \rightarrow Y$, inducing $G$-equivariant maps 
$f_\pm: {\mathcal S}_{\pm} \rightarrow X$.

\begin{lem}
\label{limite=attractif}
For any $\gamma$ in $\Gamma$ and any attractive fixed point $x$ of $\gamma$ in 
${\mathcal S}_\pm$, the image of $x$ by $f_\pm$ is an attractive 
(resp, repulsive) fixed point of $\rho(\gamma)$ in ${X}$.
\end{lem}

\preu
Quite straightforward. See proposition $3.2$ of
\cite{labourieanosov}. \fin

\begin{prop}
\label{Xattire}
The representation $\rho$ is discrete and faithful. For every non-trivial $\gamma$ in $\Gamma$,
the image $\rho(\gamma)$ is loxodromic, i.e., admits three 
eigenvalues with distinct norms.
\end{prop}

\preu
Except for the discreteness, the proposition follows immediately 
from Lemma~\ref{limite=attractif}, the fact that any non-trivial element 
of $\Gamma$ admits an attractive fixed in ${\mathcal S}_+$, and the 
fact that loxodromic elements of $\mbox{SL}(V)$ are precisely elements 
admitting an attractive fixed point in $X$. The discreteness follows 
by classical arguments. See \cite{labourieanosov} for more details.\fin

Actually, loxodromic elements of $G$ have one and only one attractive fixed 
point in $X$. Since attractive fixed points of elements of $\gamma$ 
are dense in ${\mathcal S}_\pm$, we obtain:

\begin{cor}
\label{lambda=lambda}
The maps $f_\pm$ have the same image, which is the closure of the set of attractive fixed points of elements of $\rho(\Gamma)$ in $X$.\fin
\end{cor}

\begin{defin}
\label{deflambda}
The common image $f_+({\mathcal S}_+) = f_-({\mathcal S}_-)$ is denoted by $\Lambda$, and called the limit curve.
\end{defin}

Recall that $\Phi^t$ is a double covering of the geodesic flow: 
there is a double covering between the associated leaf spaces 
${\mathcal S}_\pm \rightarrow \overline{\mathcal S}_\pm$. Let 
$\bar{\tau}$ be the Galois automorphism of this double covering. 
For any non-trivial element $\gamma$ of $\Gamma$ and any attractive 
fixed point $x$ of $\gamma$ in ${\mathcal S}_\pm$, the image 
$\bar{\tau}(x)$ is an attractive fixed point of $\gamma$. 
By uniqueness of attractive fixed points in $X$, and according 
to Lemma~\ref{limite=attractif}, we have $f_\pm(x) = f_\pm(\bar{\tau}(x))$. 
By density of attractive fixed points in ${\mathcal S}_\pm$, we obtain 
$f_\pm = f_\pm \circ \bar{\tau}$. Hence:

\begin{cor}
The maps $f_\pm$ induce maps 
$\bar{f}_\pm: \overline{\mathcal S}_\pm \rightarrow X$.\fin
\end{cor}

According to \S~\ref{bifoliated}, the orbit space $\overline{Q}$ can be 
identified with the complement $\overline{\mathcal Q}$ 
in $\overline{\mathcal S}_{+} \times \overline{\mathcal S}_-$ 
of the graph of a homeomorphism 
$\alpha: \overline{\mathcal S}_+ \rightarrow \overline{\mathcal S}_-$. 
The maps $\bar{f}_\pm$ induce a map 
$\bar{F}: \overline{\mathcal Q} \rightarrow {X} \times {X}$.

\begin{lem}
\label{imageinY}
The image of $\bar{F}$ is contained in the image $\mathcal Y$ of $\Pi$.
\end{lem}

\preu
Indeed, the maps $\bar{f}_\pm$ arise from a map $f: H \rightarrow Y$.\fin

\begin{lem}
\label{Lss}
We have the identity: $\bar{f}_+ = \bar{f}_- \circ \alpha$.
\end{lem}

\preu
Let $x$ in $\overline{\mathcal S}_+$ and $y$ in 
$\overline{\mathcal S}_- \setminus \{ \alpha(x) \}$. Then $(x,y)$ belongs 
to the open set 
$\overline{\mathcal Q} \subset \overline{\mathcal S}_+ \times \overline{\mathcal S}_-$. The image of $(x,y)$ by $\bar{F}$ 
belongs to $\mathcal Y$ (lemma~\ref{imageinY}). According to Lemma~\ref{defY}, 
we have $\bar{f}_-(y) \neq \bar{f}_+(x)$. Hence, 
$\bar{f}_+(x)$ belongs to 
$\Lambda \setminus \bar{f}_-(\overline{\mathcal S}_- \setminus \{ \alpha(x) \})$. But 
$\Lambda \setminus \bar{f}_-(\overline{\mathcal S}_- \setminus \{ \alpha(x) \})$
 is either empty or reduced to $\{ \bar{f}_-(\alpha(x)) \}$. The lemma follows.
\fin

In the proof above, we have shown in particular that $\Lambda \setminus \bar{f}_-(\overline{\mathcal S}_- \setminus \{ \alpha(x) \})$ is not empty. Hence:

\begin{cor}
The maps $\bar{f}_+$ and $\bar{f}_-$ are injective.\fin
\end{cor}

The flag manifold $X$ is a closed subset of $P(V) \times P(V^\ast)$. Let
$\eta_\pm$ (resp. $\eta^\ast_\pm$) be the composition of $\bar{f}_\pm$ with 
the projection of $X$ on $P(V)$ (resp. $P(V^\ast)$).

\begin{lem}
\label{Lfermesimple}
The maps $\eta_\pm: \overline{\mathcal S}_\pm \rightarrow P(V)$ and the maps
$\eta_\pm^\ast: \overline{\mathcal S}_\pm \rightarrow P(V^\ast)$ are injective.
\end{lem}

\preu
We only deal with $\eta_+$; the other cases are similar. 
Let $x$, $x'$ be two elements of $\overline{\mathcal S}_+$ 
with the same image by $\eta_+$: if $[u] = \eta_+(x)$, then
$\bar{f}_+(x) = ([u], [u^\ast])$, and $\bar{f}_+(x') = ([u], [v^\ast])$
for some $[u^\ast]$, $[v^\ast]$ in $P(V^\ast)$ such that 
$\langle u \mid u^\ast \rangle = \langle u \mid v^\ast \rangle =0$.

Assume $x \neq x'$. The pair $(x,y)$, with $y = \alpha(x')$, is
an element of $\overline{\mathcal Q} \approx \overline{Q}$. Its image by $\bar{F}$ is
$(\bar{f}_+(x), \bar{f}_-(\alpha(x')))$, which, according to Lemma~~\ref{Lss}, 
is equal to $(\bar{f}_+(x), \bar{f}_+(x'))$. On one hand, this pair must belong to the
image $\mathcal Y$ of $\Pi$. On the other hand, it has the form 
$(([u], [u^\ast]), ([u], [v^\ast]))$. From the description of the image of 
$\mathcal Y \subset {X} \times {X}$ (Lemma~\ref{defY}), we obtain a 
contradiction.

Hence, $x =x'$. The lemma is proved.\fin

\begin{defin}
Let $L$ be the image of $\eta_+$: this is the image of $\eta_-$ too.
Let $L^\ast$ be the common image of $\eta^\ast_+$ and $\eta^\ast_-$.
\end{defin}

According to Lemma~\ref{Lfermesimple}, $L$ and $L^\ast$ are closed
simple curves.

\begin{lem}
\label{Lambda=LL}
The limit curve $\Lambda$ is the set of pairs $([v], [v^\ast])$ where
$[v]$ belongs to $L$ and $[v^\ast]$ belongs to $L^\ast$.
\end{lem}

\preu
One of the inclusion is obvious. Conversely, let $([v], [v^\ast])$ be an element
of $X$ with $[v] \in L$, $[v^\ast] \in L^\ast$. Let $(x,y)$ be the element 
of $\overline{\mathcal S}_+ \times \overline{\mathcal S}_-$ satisfying:

\[ \eta_+(x) = [v], \;\;\;\;\; \eta_-^\ast(y) = [v^\ast] \]

Since $\langle v \mid v^\ast \rangle = 0$, the pair $(\bar{f}_+(x), \bar{f}_-(y))$
does not belong to ${\mathcal Y}$.
Therefore, $(x,y)$ cannot belong to $\overline{\mathcal Q}$. We have 
$y = \alpha(x)$.

According to Lemma~\ref{Lss}, we have $\bar{f}_+(x) = \bar{f}_-(y)$. Hence, the
$P(V^\ast)$-component of $\bar{f}_+(x)$ is $\eta_-^\ast(y) = [v^\ast]$. By construction,
the $P(V)$-component of $\bar{f}_{+}(x)$ is $[v]$. Therefore, the pair 
$([v], [v^\ast])  = \bar{f}_+(x)$ belongs to $\Lambda$. The lemma follows.\fin

By Lemma~\ref{Lambda=LL} and Lemma~\ref{Lfermesimple}:

\begin{cor}
\label{lemflag1}
Every projective line in $P(V)$ belonging to $L^\ast$ (i.e., of the form $[(u^\ast)^\perp]$ with
$[u^\ast]$ in $L^\ast$) intersects $L$ in one and only one point.\fin
\end{cor}

The dual statement, with $L$ and $L^\ast$ exchanged, is of course true. Hence, the corollary
above can be complemented by:

\begin{cor}
\label{lemflag2}
For every point $[u]$ in $L$, there is one and only one projective line of $P(V)$ belonging to
$L^\ast$ and containing $[u]$.\fin
\end{cor}

\begin{cor}
\label{Lattire}
The curve $L$ (respectively $L^\ast$) is the closure of the set of attractive fixed points in
$P(V)$ (respectively $P(V^\ast)$) of elements of $\rho(\Gamma)$ (respectively
$\rho^\ast(\Gamma)$).\fin
\end{cor}

\rque
Of course, in the canonical case, i.e., when $\rho$ is the restriction of $\rho_0$ to 
$\Gamma \subset H$, the curves $L$, $L^\ast$ are the projective lines $L_0$, $L_0^\ast$.
In this case, we defined Klein bottles $T_0$, $T^\ast_0$ and an open domain $\Omega_0$ 
(see \S~\ref{exempleflagmanifold}). These constructions extend to the general case in 
the following way: define $T$
(respectively $T^\ast$) to be the set of flags $([u], [u^\ast])$ with $[u] \in L$ (respectively $[u^\ast] \in L^\ast$), and let $\Omega$ be the complement in
$X$ of the union $T \cup T^\ast$. 

According to Lemma~\ref{Lambda=LL}, the limit curve $\Lambda$ is the intersection 
$T \cap T^\ast$, and the complements of $\Lambda$ in $T$, $T^\ast$ are denoted $A$, $A^\ast$.
\erque

\rque
\label{fauxlorentz}
All the results of these sections apply to any $(G,Y)$-Anosov representation and, in particular,
to quasi-Fuchsian representations. In this case, $L$ is the boundary of a $\rho(\Gamma)$-invariant convex domain $C$ in $P(V)$, and $L^\ast$ is the boundary of the dual convex
$C^\ast$: elements in $L^\ast$ are projective lines in $P(V)$ tangent to $L$. The sets $T$, $T^\ast$ are (topological) tori, and $A$, $A^\ast$ are annuli.

The domain $\Omega$ in this case has $3$ connected components:

-- one component is the set of flags $([v], [v^\ast])$ with $[v] \in C$: this component is canonically identified with the projectivized tangent bundle of $C$,

-- another component is the set of flags $([v], [v^\ast])$ with 
$[v^\ast] \in C^\ast$: it is canonically identified with the projectivized tangent bundle of $C^\ast$,

-- Finally, there is a third component, consisting of the flags $([v], [v^\ast])$ 
with $[v] \notin C$, $[v^\ast] \notin C^\ast$.

The last component, in some way, has a lorentzian flavor. Indeed, when $\rho$ is Fuchsian,
i.e., when $C$ and $C^\ast$ are ellipses, this last component is canonically identified with
the projectivized bundles of timelike vectors of de Sitter space. 
\erque

\section{Special deformations}

In this section, we fix the embedding $\Gamma \subset H$, i.e., the canonical morphism 
$\rho_0: \Gamma \rightarrow G$. The projection of $\Gamma$ in $\overline{H}$ 
is injective; we still denote by $\Gamma$ the image of this projection. 
The quotient $\Gamma\backslash\overline{H}$ is the surface $\Sigma$. 

\begin{defin}For any $\gamma$ in $\Gamma$ let $r(\gamma)$ be the spectral 
radius of $\gamma \in H$.
\end{defin}

Up to the sign, the eigenvalues of $\gamma \in H$ are
$r(\gamma)$, $r(\gamma)^{-1}$; when $\gamma$ is non-trivial, we have $r(\gamma) >1$.

\subsection{$\Lambda$-preserving deformations}
\label{flagdesargue}

\subsubsection{A flag version of the geodesic flow}
\label{varphhi}
Consider\footnote{The factor $1/3$ arises from the fact that 
this flow truely lies in $\mbox{PGL}(V)$:
multiplying every coefficient by the inverse of the middle diagonal coefficient provides 
a more elegant expression{\ldots} } the following $1$-parameter subgroup of $G$:

\[ \varphi^t = \left(\begin{array}{ccc}
       e^{t/3} & 0 & 0 \\
       0 & e^{-2t/3} & 0 \\
       0 & 0 & e^{t/3}
       \end{array}\right) \]

It commutes with every $\rho_0(h)$. Hence, it preserves the open 
$\rho_0(H)$-orbit $\Omega_0$ in $X$ and induces a flow on the quotient manifold 
$\rho_0(\Gamma)\backslash\Omega_0$. The projection of
$\Omega_0$ on $P(V)$ is contained in the complement of $L_0$, an affine plane. 
Consider the coordinate system $(u,v)$ on this plane such that the
coordinates of $[ue_1 + e_2 + ve_3]$ are $(u,v)$.

Actually, the projection of $\Omega_0$ in $P(V)$ is the complement 
in this plane of the point $[e_2]$, i.e., the complement of the origin $(u,v) = (0,0)$. 
The induced action of $\varphi^t$ on this projection is the homothety
of factor $e^t$ fixing $[e_2] = (0,0)$. 

On the other hand, if $h$ is an element of $H \approx \mbox{SL}(2, {\mathbb R})$ of the form:

\[ \left(\begin{array}{cc}
         a & b \\
	 c & d
	 \end{array}\right) \]

then $f_+(h)$ projects to $[\rho_0(h)e_1]$ in $P(V)$, with $(u,v)$-coordinates $(a,c)$. 
Hence, for any $t$, the $(u,v)$-coordinates of $f_+(ha^t)$ are $(e^ta, e^tc)$.
Since the action of a transformation of $X$ is characterized by its projective action 
on $P(V)$, we see that, via the identification $\Omega_0 \approx H$ we have selected, 
\emph{the action of $\varphi^t$ on $\Omega_0$ coincides with the right action of $a^t$ on $H$.\/} 
Hence, \emph{the flow on $\rho_0(\Gamma)\backslash\Omega_0$ induced by $\varphi^t$ is 
conjugate
to the flow $(M, \Phi^t)$.\/}

\subsubsection{Linear deformations}

Let $u: \Gamma \rightarrow \mathbb R$ be any morphism. The canonical morphism
can be deformed to a new morphism, called the \emph{$u$-deformation:\/}

\[ \rho_u(\gamma) = \varphi^{u(\gamma)} \circ \rho_{0}(\gamma) \]

The morphism $u$ is an element of $H^{1}(\Gamma, {\mathbb R})$. On this cohomology space, 
with $\rho_0: \Gamma \rightarrow H$ fixed, we can define the \emph{stable norm\/} 
(cf. \cite{norm}) as
follows: for any hyperbolic element $\gamma$ of $\Gamma$,
let $t(\gamma)$ be the double of the logarithm of $r(\gamma)$
(this is the length of the closed geodesic associated to $\Gamma$
in the quotient of the Poincar\'e disc by $\Gamma$).
For any element $\hat{\gamma}$ of $H_{1}(\Gamma, {\mathbb Z})$,
and for any positive integer $n$, let $t_{n}(\hat{\gamma})$ be the
infimum of the values $\frac{t(\gamma)}{n}$ where $\gamma$ describes
all the elements of $\Gamma$ representing $n\hat{\gamma}$.
The limit of $t_{n}(\hat{\gamma})$ exists; it is the {\em stable
norm of $\hat{\gamma}$ in $H_{1}(\Gamma, {\mathbb Z})$.\/} 
This norm is extended in a unique way on all 
$H_{1}(\Gamma, {\mathbb R})$; the dual of it is the
{\em stable norm of $H^{1}(\Gamma, {\mathbb R})$.\/}
The stable norm of $u$ in  $H^{1}(\Gamma, {\mathbb R})$ is denoted $\vert u \vert_s$.

\begin{thm}
\label{flaganosov}
The representation $\rho_u$ is $(G,Y)$-Anosov if and only if $\vert u \vert_s < 1/2$. 
\end{thm}

\preu
Assume that $\rho_u$ is $(G,Y)$-Anosov. The invariant curve $\Lambda$ must be the closure 
of the union of attractive fixed points of $\rho_u(\Gamma)$. In particular, it contains the closure 
of the attractive fixed points of the commutator subgroup $\rho_u([\Gamma, \Gamma])$.
But $\rho_0$ and $\rho_u$ coincide on $[\Gamma, \Gamma]$, and attractive fixed points 
of elements of $[\Gamma, \Gamma]$ are dense in $\overline{\mathcal S}_\pm \approx {\mathbb R}P^1$; 
hence, $\Lambda = \Lambda_0$.

But it is easy to see that if $\vert u \vert_s > 1/2$, there is an element $\gamma$ of $\Gamma$ 
such that the attractive fixed point of $\rho_u(\gamma)$ in $P(V)$ is $[e_2]$. This attractive 
fixed point does not belong to $L$, which contradicts Corollary~\ref{Lattire}.
Therefore, we must have $\vert u \vert_s \leq 1/2$. Since the $(G,Y)$-Anosov property is open in $\mbox{Rep}(\Gamma, G)$, the inequality is strict: $\vert u \vert_s < 1/2$.

The inverse statement, i.e., the fact that $\rho_u$ is $(G,Y)$-Anosov if $\vert u \vert_s < 1/2$, 
is proved in Appendix A.
\fin

\rque
As the proof above shows, the limit curve $\Lambda$ of $\rho_u$, when 
$\vert u \vert_s < 1/2$, is the limit curve $\Lambda_0$ of $\rho_0$. 
It does not depend on the inclusion $\Gamma \subset H$.
\erque

\subsection{Deformations with $L_0$ remaining constant}
\label{flagoldman}

Consider morphisms $\rho: {\Gamma} \rightarrow G$ of the form:

\[ \rho(\gamma) = \left(\begin{array}{ccc}
                  e^{u(\gamma)/3}a(\gamma) & 0 & e^{u(\gamma)/3}b(\gamma) \\
		  \mu(\gamma) & e^{-2u(\gamma)/3} & \nu(\gamma) \\
		  e^{u(\gamma)/3}c(\gamma) & 0 & e^{u(\gamma)/3} d(\gamma)
\end{array}\right) \]
where:

--  $\rho_\lambda: \gamma \mapsto \left(\begin{array}{cc}
                    a(\gamma) & b(\gamma) \\
		    c(\gamma) & d(\gamma) 
		    \end{array}\right) $ is a Fuchsian representation taking value in $\mbox{SL}(2,{\mathbb R})$, i.e., is injective, with discrete image.

--  $u: {\Gamma} \rightarrow {\mathbb R}$ is a morphism.

Such a representation is called a \emph{radial representation.\/}

When $u$ has stable norm (relatively to $\rho_\lambda$) strictly less than $1/2$, 
$\rho$ is called a \emph{hyperbolic representation.\/} In this case, $[e_2]$ is a fixed point of saddle type of every non-trivial $\rho(\gamma)$.
In \cite{barflag}, the action of $\Gamma$ on $P(V)$ induced by such a representation is called a \emph{hyperbolic action.\/}

One of the main results of \cite{barflag} is:

\begin{thm}[Theorem A of \cite{barflag}]
Let $\rho$ be a hyperbolic representation.
The action of $\rho(\Gamma)$ on $P(V)$ preserves a continuous closed simple curve $L$. This curve, if  Lipschitz, is a projective line.
\end{thm}

Actually, the first part of this Theorem is a corollary of the following:

\begin{lem}
\label{chouia}
Hyperbolic representations are $(G,Y)$-Anosov. 
\end{lem}

\preu
If the maps $\mu$ and $\nu$ are trivial, it follows from Theorem~\ref{flaganosov}. 
Then, observe that the composition of a hyperbolic representation by the conjugacy 
in $G$ by $\varphi^t$ remains hyperbolic radial, with the same $\rho_\lambda$ and 
the same $u$, but with coefficients $\mu$, $\nu$ multiplied by $e^{-2t}$. 
Since $\rho_u$ is $(G,Y)$-Anosov, and by stability of $(G,Y)$-Anosov representations, 
the conjugated representation for big $t$ is $(G,Y)$-Anosov. The lemma follows since 
conjugacy in $G$ does not affect the $(G,Y)$-Anosov property.\fin

\rque
The invariant curve $L^\ast$ in $P(V^\ast)$ is obviously $L_0^\ast = [e_2^\perp]$.
\erque

\rque
The analog of Theorem~\ref{flaganosov} remains true: a radial representation is $(G,Y)$-Anosov if and only if $\vert u \vert_s < 1/2$. The proof is similar; we leave it to the reader.
\erque

\begin{prop}
\label{pointfixe?}
Let $\rho: \Gamma \rightarrow G$ be a $(G,Y)$-Anosov representation. Assume that a finite index subgroup of $\rho(\Gamma)$ preserves a proper subspace of $V$. Then, up to conjugacy in $G$, $\rho$ or $\rho^\ast$ is a hyperbolic representation.
\end{prop}

\preu
A proper subspace of $V$ is a line or a $2$-plane: replacing $\rho$ by $\rho^\ast$ if necessary, we can assume that the subspace preserved by a finite index subgroup of $\rho(\Gamma)$ is a line. After conjugacy in $G$, we can assume moreover that this invariant line is spanned by $e_2$. Let $\Gamma' \subset \Gamma$ be the finite index subgroup such that $\rho(\Gamma')$ fixes $[e_2]$. The restriction $\rho'$ of $\rho$ to $\Gamma'$ is still $(G,Y)$-Anosov, and its limit curve $\Lambda'$ is the limit curve $\Lambda$ of $\rho$.

According to Proposition~\ref{Xattire}, $\rho'$ is faithful, with discrete image. It follows that $\rho'$ is a radial representation, described as above by a Fuchsian representation $\rho'_\lambda: \Gamma' \rightarrow H$ and a morphism $u': \Gamma \rightarrow \mathbb R$. This morphism is trivial
on $[\Gamma', \Gamma']$. It follows that $[e_2]$ is a saddle fixed point of 
$\rho'(\gamma)$ for every element $\gamma$ of $[\Gamma', \Gamma']$. 
Hence, the attractive fixed point of $(\rho')^\ast(\gamma)$ in $P(V^\ast)$ 
belongs to $[e_2^\perp]$. The argument used in the proof of Theorem~\ref{flaganosov} and 
the identity $\Lambda = \Lambda'$ imply here that the invariant curve $L_0^\ast$ is $[e_2^\perp]$.
In particular, it follows that $[e_2]$ is a saddle fixed point of $\rho(\gamma)$ 
for every $\gamma$ in $\Gamma$, not only in $[\Gamma', \Gamma']$. 
Therefore, $\rho$ is a linear $u$-deformation of a canonical representation. 
As in the proof of \ref{flaganosov}, we infer $\vert u \vert_s \leq 1/2$, and, finally, 
$\vert u \vert_s < 1/2$ (by stability of Anosov representations).\fin

\section{Properness of the action on $\Omega$}
\label{propre}

Let $\rho: {\Gamma} \rightarrow G$ be a $(G,Y)$-Anosov representation. Recall that $\rho(\Gamma)$ preserves the open domain $\Omega$, which is the set of flags $([u], [u^\ast])$ with $[u] \notin \L$, $[u^\ast] \notin L^\ast$.

\begin{thm}
\label{actiopropre}
For any $(G,Y)$-Anosov representation $\rho: {\Gamma} \rightarrow G$, the action of $\rho(\Gamma)$ action of ${\Gamma}$ on the associated domain $\Omega$ is free and properly discontinuous.
\end{thm}

\rque
\label{casreduc}
This Theorem, when $\rho$ is a linear deformation $\rho_u$, is a reformulation of Th\'eor\`eme $3.4$ of \cite{salein}. More generally, Theorem~\ref{actiopropre} is proved in \cite{barflag}
(Proposition $4.19$ of \cite{barflag})
in the case of hyperbolic representations.

Therefore, in this section, we don't consider this particular case. By  
Proposition~\ref{pointfixe?}, it means that we assume that no proper subspace 
of $V$ is $\rho(\Gamma)$-invariant, i.e., that $\rho$ is \emph{strongly irreducible.\/}

In Appendix B, we show how the proof developed here can be adapted to the reducible case, achieving a complete proof of \tm~\ref{actiopropre}.
\erque

All this section is devoted to the proof of \tm~\ref{actiopropre}. Let $\rho$ be a 
$(G,Y)$-Anosov representation. According to Lemma~\ref{Xattire}, every $\rho(\gamma)$ 
is loxodromic: its fixed points in $X$ all have the form $([u], [u^\ast])$ where either 
$[u]$ belongs to $L$, or $[u^\ast]$ belongs to $L^\ast$. None of them belongs to $\Omega$: 
the action
of $\rho(\Gamma)$ is free.

Our task is to prove that the action is proper. Since $\rho(\Gamma)$ is discrete, 
it amounts to proving that there is no sequence $(\gamma_{n})_{(n \in {\mathbb N})}$ 
in $\Gamma$ for which there is a sequence of flags $([u_n], [u^\ast_n])_{(n \in {\mathbb N})}$ 
in $\Omega$ satisfying:

\begin{enumerate}

\item the sequence $\rho(\gamma_n)_{(n \in {\mathbb N})}$ escapes from any compact of $G$,

\item the sequence $([u_n], [u^\ast_n])_{(n \in {\mathbb N})}$ converge to a point $([\bar{u}],
[\bar{u}^\ast])$ in $\Omega$,

\item the sequence 
$\rho(\gamma_n)([u_n], [u^\ast_n])_{(n \in {\mathbb N})} = ([v_n], [v^\ast_n])_{(n \in {\mathbb N})}$ 
converges to a flag $([\bar{v}], [\bar{v}^\ast])$ in $\Omega$.

\end{enumerate}

We argue by contradiction, assuming the existence of such a sequence. We denote 
$g_n = \rho(\gamma_n)$.

Remember that we equip $V$, $V^\ast$ with the euclidean metrics $N$, $N^\ast$ for which 
$([e_1], [e_2], [e_3])$ and $([e_1^\ast], [e_2^\ast], [e_3^\ast])$ are an orthonormal basis. 
It induces a norm on $\mbox{GL}(V)$: the operator norm $\Vert$. 
Let ${\mathcal A} \subset \mbox{GL}(V)$ be the unit sphere of this norm.

Consider the Cartan decompositions of $g_n$, $g_n^\ast$ according to the canonical basis:

\[ g_n = k_n \left(\begin{array}{ccc}
             \lambda_n & 0 & 0 \\
	     0 & \mu_n & 0 \\
	     0 & 0 & \nu_n
	     \end{array}\right) l^{-1}_{n} \]

\[g^\ast_n = k_n \left(\begin{array}{ccc}
             \lambda_n^{-1} & 0 & 0 \\
	     0 & \mu_n^{-1} & 0 \\
	     0 & 0 & \nu_n^{-1}
	     \end{array}\right) l^{-1}_{n} \]
	     
where $k_n$, $l_n$ are isometries of $N$, and  
$\lambda_n \geq \mu_n \geq \nu_n$, with $\lambda_n\mu_n\nu_n = 1$.

The quotients $\bar{g}_n = \frac{g_n}{\Vert g_n \Vert}$ and 
$\bar{g}_n^\ast = \frac{g_n^\ast}{\Vert g^\ast_n \Vert}$ belong to $\mathcal A$. 
By compactness of this unit sphere, passing to a subsequence if necessary, we can 
assume that $\bar{g}_n$ and $\bar{g}_n^\ast$ converge to elements $\bar{g}$, $\bar{g}^{\ast}$,
repectively,  of $\mathcal A$.
Let $I$ be the image of $\bar{g}$, $I_\ast$ the image of $\bar{g}^\ast$, $K$ the kernel of $\bar{g}$ and $K_\ast$ the kernel of $\bar{g}^{\ast}$. By convention, $[I]$,
$[K]$, $[I_\ast]$ and $[K_\ast]$ are their projections in $P(V)$, $P(V^\ast)$.
Let $I^{\perp}, K^{\perp} \subset V^\ast$ and 
$K^{\perp}_{\ast}, I^{\perp}_\ast \subset V$ be the dual subspaces, and $[I^\perp]$, $[K^\perp]$, {\ldots} their respective projectivizations.

The $g_n$ all have determinant $1$. Hence, the item $(1)$ above implies that, 
after taking a subsequence, the norms $\Vert g_n \Vert$ and $\Vert g_{n}^{\ast} \Vert$ tend to 
$+\infty$. If 
$\bar{g}^t$ denotes the transpose matrix of $\bar{g}$, we have:

\[ \bar{g}^t \circ \bar{g}^\ast = 0 = \bar{g}^\ast \circ \bar{g}^t \]

Therefore:

\[ I_\ast \subset I^\perp, \;\;\;\;\;\; K^\perp \subset K_\ast \]

Since the actions of $g_n$ and $\bar{g}_n$ on $P(V)$ coincide, as for the actions of $g^\ast_n$ 
and $\bar{g}^\ast_n$ on $P(V^\ast)$, we obtain:

\begin{lem}
\label{caconverge}
The sequence $(\bar{g}_{n})_{(n \in {\mathbb N})}$ converge uniformly on compact subsets of
$P(V) \setminus [K]$ to the restriction of $\bar{g}$ to $P(V) \setminus [K]$. 
Similarly, the sequence $(\bar{g}^\ast_{n})_{(n \in {\mathbb N})}$ converge uniformly on 
compact subsets of $P(V^\ast) \setminus [K_\ast]$ to the restriction of $\bar{g}$.\fin
\end{lem}

\rque
The image of the restriction of $\bar{g}$ to $P(V) \setminus [K]$ is $[I]$, and the image of
the restriction of $\bar{g}^\ast$ to $P(V^\ast) \setminus [K_\ast]$ is $[I_\ast]$.
\erque

Actually, since the norm of $g_n$ tends to $+\infty$, we see that, for 
$\bar{k}$ and $\bar{l}$ denoting the limits of $\bar{k}_{n}$ and 
$\bar{l}_{n}$, the matrices $\bar{g}$ and $\bar{g}^{\ast}$ have the expressions: 

\[ \bar{g} = \bar{k} \left(\begin{array}{ccc}
             \bar{\lambda} & 0 & 0 \\
	     0 & \bar{\mu} & 0 \\
	     0 & 0 & 0
	     \end{array}\right) \bar{l}^{-1} \]

\[\bar{g}^\ast = \bar{k} \left(\begin{array}{ccc}
             0 & 0 & 0 \\
	     0 & \bar{\mu}' & 0 \\
	     0 & 0 & \bar{\nu}
	     \end{array}\right) \bar{l}^{-1} \]

where $\bar{\lambda}$, $\bar{\nu}$ are positive, and $\bar{\mu}$, $\bar{\mu}'$ are non-negative.

We will also need to exchange the roles of $([\bar{u}],
[\bar{u}^\ast])$ and $([\bar{v}], [\bar{v}^\ast])$; this amounts to replacing the 
$\gamma_{n}$ by their inverses, i.e., to replacing $g_n$ and 
$g^\ast_n$. Let $h_n$, $h^\ast_n$ be these inverses: 

\[ h_n = l_n \left(\begin{array}{ccc}
             \lambda_n^{-1} & 0 & 0 \\
	     0 & \mu_n^{-1} & 0 \\
	     0 & 0 & \nu_n^{-1}
	     \end{array}\right) k^{-1}_{n} \]

\[h^\ast_n = l_n \left(\begin{array}{ccc}
             \lambda_n & 0 & 0 \\
	     0 & \mu_n & 0 \\
	     0 & 0 & \nu_n
	     \end{array}\right) k^{-1}_{n} \]

Their limits for $n \rightarrow +\infty$ are:

\[ \bar{h} = \bar{k} \left(\begin{array}{ccc}
             0 & 0 & 0 \\
	     0 & \bar{\mu}' & 0 \\
	     0 & 0 & \bar{\nu}'
	     \end{array}\right) \bar{l}^{-1} \]

\[\bar{h}^\ast = \bar{k} \left(\begin{array}{ccc}
             \bar{\lambda} & 0 & 0 \\
	     0 & \bar{\mu} & 0 \\
	     0 & 0 & 0
	     \end{array}\right) \bar{l}^{-1} \]

\begin{lem}
\label{muestnul}
The coefficients $\bar{\mu}$ and $\bar{\mu}'$ are both zero.
\end{lem}

\preu
Assume that one of them, let's say $\bar{\mu}$, is non-zero. Then, $K$ is
the line spanned by $\bar{l}(e_3)$. The inclusion
$K^\perp \subset K_\ast$ shows that $K_\ast$ has dimension $2$. Hence, 
$I_\ast$ is $1$-dimensional. Therefore, the inclusions are all identities:

\[ I_\ast = I^\perp, \;\;\;\;\;\; K^\perp = K_\ast \]

\emph{Fact $1$: $[K]$ belongs to the invariant curve $L$.\/}
If not, the $g_n$-invariant curve $L$ is contained in $P(V) \setminus
[K]$. By  Lemma~\ref{caconverge}, it implies that $L$ is equal to the projective line $[I]$, 
which is therefore a $\rho(\Gamma)$-invariant projective line.
This contradicts the strong irreducibility of $\rho$.

\emph{Fact $1'$: $[I_\ast]$ belongs to $L^\ast$.\/}
Apply the argument used for the proof of Fact $1$ to the inverse sequence
$(h^\ast = (g_n^\ast)^{-1})_{(n \in {\mathbb N})}$. We obtain that if $[I^\ast]$ 
does not belong to $L^\ast$, then $L^\ast$ is the projective line $[K_\ast]$, contradicting 
once more the irreducibility of $\rho$.

\emph{Fact $2$: The invariant curve $L$ is a projective line containing $[K]$.\/}
Indeed, according to Fact $1'$, the projective line $[I]$ belongs to $L^\ast$. According 
to Corollary~\ref{lemflag1}, the intersection $[I] \cap L$ is a single point $[p]$. 
On the other hand, the fibers of the restriction of $\bar{g}$ to $P(V) \setminus [K]$ 
are the projective lines containing $[K]$, with the point $[K]$ excluded. 
By Lemma~\ref{caconverge}, since $[I] \cap L$ is a single point, $L \setminus [K]$ 
is contained in one and only one of these fibers. Fact $2$ follows.

The lemma follows by irreducibility of $\rho$.\fin

According to Lemma~\ref{muestnul}, $[K]$ is a projective line, and $[I]$ is a single point. Similarly, $[K_\ast]$ is a projective line, whereas 
$[I_\ast]$ is a single point.

\begin{lem}
\label{IIast}
The points $[I]$, $[I_\ast]$ belong to $L$, $L^\ast$, respectively.
\end{lem}

\preu
Assume that $[I]$ does not belong to $L$. Then, by Lemma~\ref{caconverge}, 
$L$ is contained in the projective line $[K]$, which is therefore 
$\rho(\Gamma)$-invariant, a contradiction.
The proof of the dual statement $[I_\ast] \in L^\ast$ is similar.\fin

\begin{lem}
\label{KKast}
The points $[K^\perp]$ and $[K_\ast^\perp]$ belong respectively to $L^\ast$, $L$.
\end{lem}

\preu
 According to Lemma~\ref{muestnul},  the limit in ${\mathcal A}$ of the $\bar{h}_{n}$ (i.e., the $h_n$ divided by their norms) is:

\[\bar{h} = \bar{l} \left(\begin{array}{ccc}
             0 & 0 & 0 \\
	     0 & 0 & 0 \\
	     0 & 0 & \bar{\nu}
	     \end{array}\right) \bar{k}^{-1} \]

Hence, the image of $\bar{h}$ is spanned by $\bar{l}e_3$. On the 
other hand, the kernel  $K_\ast$ of $\bar{g}^\ast$ is the $2$-plane 
generated by $\bar{l}(e^\ast_1)$ and $\bar{l}(e^\ast_2)$. 
Therefore, $K_\ast^\perp$ is the line spanned by $\bar{l}e_3$, i.e., 
the image of $\bar{h}$. Applying Lemma~\ref{IIast} to $\bar{h}$, 
we get that $[K^\perp_\ast]$ belongs to $L$. 

The same argument, applied to $\bar{h}^\ast$, implies $[K^\perp] \in L^\ast$ \fin

{\bf \mbox{End of the Proof of Theorem~\ref{actiopropre} in the strongly irreducible case: }}

Since $([\bar{v}], [\bar{v}^\ast])$ belongs to $\Omega$, $[\bar{v}]$ does not belong to $L$.
According to Lemma~\ref{IIast}, $\bar{v} \neq [I]$. It follows, with Lemma~\ref{caconverge},that $[\bar{u}]$ belongs to $[K]$. Similarly, we have $[\bar{u}^\ast] \in [K_\ast]$. Hence:

-- the projective line $[K_\ast]$ in $P(V^\ast)$ contains $[\bar{u}^\ast]$ and $[K^\perp]$ (since $K^\perp \subset K_\ast$),

-- the projective line $[\bar{u}^{\perp}]$ in $P(V^\ast)$ contains $[\bar{u}^\ast]$ (since $ \langle \bar{u} \mid \bar{u}^\ast \rangle =
0$) and $[K^\perp]$ (since $[\bar{u}] \in [K]$) too. 

But, according to Lemma~\ref{KKast}, since $[\bar{u}^\ast]$ does not belong to  $L^\ast$, 
we have $\bar{u}^\ast \neq [K^\perp]$. Hence the projective lines 
$[K_\ast]$ and $[\bar{u}^\perp]$ share two distinct points: 
they are equal. In other words, $[\bar{u}]$ is equal to $[K_\ast^\perp]$: 
according to Lemma~\ref{KKast}, it belongs to $L$, a contradiction.\fin

\rque
\label{fuchspropre}
\tm~\ref{actiopropre} is true for any $(G,Y)$-Anosov representations, in parti\-cular, it applies 
to representations in the Hitchin component, i.e., dividing a strictly convex domain 
$C$ of $P(V)$.
In this case, we recover the well-known properness of the action on the projectivized 
tangent bundle of $C$. But we also obtain the properness of the action on 
"timelike directions" over the "de Sitter-like" component $P(V) \setminus \overline{C}$.
\erque

\rque
We can ask how \tm~\ref{actiopropre} can be extended to any other $(G,Y)$-Anosov representations, 
for other pairs $(G,Y)$. For example, for the pair 
$(G,Y) = (\mbox{PSL}(2,{\mathbb C}), {\mathbb C}P^{1} \times {\mathbb
C}P^1 \setminus \mbox{diag})$, i.e., the case of conformal quasi-Fuchsian representations, 
\tm~\ref{actiopropre} corresponds the well-known properness of the action on the discontinuity 
domain, i.e., the complement in ${\mathbb C}P^1$ of the limit set.
\erque

\rque
In the wonderful paper \cite{pappus}, R. Schwarz considers some particular actions of 
$\mbox{PSL}(2, {\mathbb Z})$ on the flag manifold $X$. He exhibited invariant curves 
$L$, $L^\ast$ with properties completely similar to the invariant curves $L$, $L^\ast$ we have considered here\footnote{There is a minor difference: the action on $X$ considered in \cite{pappus} contains \emph{polarities,\/} i.e., projective transformations followed by the flip 
$([v], [v_\ast]) \mapsto ([v^\perp_\ast], [v^\perp])$.\/}: compare our corollaries 
\ref{lemflag1}, \ref{lemflag2} with Theorem
$3.3$ of \cite{pappus}. Our Theorem~\ref{actiopropre} echoes Theorem $4.2$ of \cite{pappus}.
Hence, it seems reasonable to qualify these representations of $\mbox{PSL}(2, {\mathbb Z})$
as $(G,Y)$-Anosov representations. In fact, they are also deformations of the representation obtained by composing  $\rho_{0}: H
\rightarrow G$ with the inclusion $\mbox{SL}(2, {\mathbb Z}) \subset \mbox{SL}(2, {\mathbb R})$. But there is a crucial difference: 
$\mbox{PSL}(2, {\mathbb Z})$ is indeed a lattice of $\mbox{PSL}(2, {\mathbb R})$, 
but not cocompact! 
It is therefore presumably possible and interesting to extend the notion of 
$(G,Y)$-Anosov representations to nonuniform lattices of $H$.
\erque

\section{Anosov flag manifolds}
\label{anosovmanifold}

Thanks to Theorem~\ref{actiopropre}, the definition~\ref{defcan} can be extended:

\begin{defin}
\label{deflaganosov}
Let $\rho: {\Gamma} \rightarrow G$ be a $(G,Y)$-Anosov representation. 
The quotient by $\rho({\Gamma})$ of the domain $\Omega$ is called an Anosov flag manifold..
\end{defin}

\subsection{Tautological foliations}
\label{tautogoldman}
The fibers of the projections $X \rightarrow P(V), P(V^\ast)$ 
are leaves of foliations by circles on $X$. They are preserved by $G$; therefore, 
they induce two $1$-dimensional foliations on any flag manifold. 
The foliation corresponding to the projection $X \rightarrow P(V)$ is called the 
\emph{first tautological foliation.\/} The other one, corresponding to the projection 
$X \rightarrow P(V^\ast)$, is called the \emph{second tautological foliation.\/} 
Observe that these foliations are both transversely real projective. 

\rque
\label{tautoFuchsian}
When $\rho$ is hyperconvex, $\Omega$ has three connected components (see remark~\ref{fauxlorentz}):

-- one is the set of flags $([v], [v^\ast])$ with $[v] \in [C]$. We denote it by $\Omega_1$.
The quotient $M_1 = \rho(\Gamma)\backslash\Omega_1$ is naturally identified with the
projectivized tangent bundle of the convex real projective surface 
$S = \rho(\Gamma)\backslash{[C]}$. The leaves of the first tautological 
foliation are the fibers of the bundle map. In particular, they are compact. 
The second tautological foliation is the foliation supported by the geodesic flow 
of the Hilbert metric on $S$, quotiented by the involution 
sending any tangent vector to its opposite: hence, it is topologically conjugate 
to the geodesic flow of any hyperbolic metric on the surface $\Sigma$ quotiented by 
the antipodal map in the fibers (see \cite{barflag} for more details).

-- another component, $\Omega_2$, which is the set of flags 
$([v], [v^\ast])$ with $[v^\ast] \in [C^\ast]$. The quotient 
$M_2 = \rho(\Gamma)\backslash\Omega_2$ is the projectivized tangent 
bundle of the dual convex real projective surface 
$S^\ast = \rho^\ast(\Gamma)\backslash{[C^\ast]}$. The leaves of the 
second tautological foliation are the fibers of the bundle map, whereas 
the first tautological foliation is doubly covered by the geodesic flow on $\Sigma$.

-- the third component, $\Omega_3$, is the ``de Sitter''-like component. It is the set 
of pairs $([v], [v^\ast])$, where $[v]$ belongs to the M\"{o}bius band ${\mathcal A}$, 
which is the complement in $P(V)$ of the closure of $[C]$, and $[v^\ast]$ belongs to 
the M\"{o}bius band ${\mathcal A}^\ast$, which is the complement in $P(V^\ast)$ of the 
closure of $[C^\ast]$. 
Observe that for any $[v]$ in ${\mathcal A}$, the intersection $[v^\perp] \cap L^\ast$ 
is the union of two points $[v^\ast_-]$, $[v_+^\ast]$ (note that 
$L^\ast = \partial [C^\ast]$). Then, $[v^\ast_\pm]$ is tangent to $L = \partial [C]$ 
at a point $[v_\pm]$. Let $p([v])$ be the intersection of $[(v^\ast)^\perp]$ and 
the projective line $[w^\ast]$ containing $[v_-]$ and $[v_+]$; it belongs to $[C]$. 
Hence, the flag $(p([v]), [w^\ast])$ belongs to $\Omega_1$. We have thus defined a 
$\rho(\Gamma)$-equivariant map from $\Omega_3$ into $\Omega_1$. It is easy to show that 
it is a homeomorphism, and that it sends the first tautological foliation of $\Omega_3$ 
onto the second tautological foliation of $\Omega_1$.

A similar treatment can be applied to the second tautological foliation on 
$\Omega_3$, leading to the following statement: \emph{the tautological 
foliations of $M_3 = \rho(\Gamma)\backslash\Omega_3$ are both topologically 
conjugate to the geodesic flow of any hyperbolic metric on the surface $\Sigma$ 
quotiented by the antipodal map in the fibers.\/}
\erque

\subsection{Canonical Goldman flag manifolds}
In \cite{barflag}, Anosov flag manifolds associated to linear $u$-deformations of 
canonical representations were defined and called 
\emph{canonical Goldman flag manifolds.\/} 

We recall that in the case of canonical Goldman flag manifolds, 
the flow $\varphi^t$ induces an Anosov flow on the associated canonical flag manifold. 
More precisely, this flow is a double covering of a \emph{Desarguian Anosov flow\/} 
as defined in \cite{ghyexotique}, \cite{barTCP}. The tautological foliations are 
then the strong stable and unstable foliations. See \cite{barflag} for more details.

For our purpose here, it is more suitable to drop the identification 
$G = \mbox{PGL}(V) \approx \mbox{SL}(V)$, and to consider 
$u$-deformations as morphisms ${\lambda}: \Gamma \rightarrow \mbox{GL}_2$, 
where $\mbox{GL}_2$ is the group of invertible $2 \times 2$ matrices, 
identified with the stabilizer in $G$ of $[(L_0^\ast)^\perp] = [e_2]$ and 
$L_0 = [(e_2^\ast)^\perp]$:

\[  \rho(\gamma) = \left(\begin{array}{ccc}
       \lambda(\gamma) & \begin{array}{c}
                        0 \\
                        0 \end{array} \\
       \begin{array}{cc}
       0 & 0 \end{array} & 1 
       \end{array}\right) \]

The group $\mbox{GL}_2$ preserves the annulus
$A_{can} = P(V) \setminus (L_0 \cup [e_2])$ (there is also a dual action on the annulus $P(V^\ast) \setminus (L_0^\ast \cup [e_2^\ast])$. The minimality and unique ergodicity of horocyclic flows (\cite{ergodic})  imply that the ${\lambda}(\Gamma)$-action on $A_{can}$ is minimal and uniquely ergodic: up to a positive factor, there is one and only Borel measure on $A_{can}$ preserved by ${\lambda}(\Gamma)$.

\begin{prop}
\label{rigid}
Two inclusions ${\lambda}_1$, ${\lambda}_2$ from a surface group $\Gamma$ 
into $\mbox{GL}_2$ induce conjugate actions on the annulus $A_{can}$ if and only 
if ${\lambda}_1$, ${\lambda}_2$ are conjugate in $\mbox{GL}_2$.
\end{prop}

\preu
This is a folkloric fact, but we don't know any appropriate reference.

One of the implication is clear. Let's prove the inverse statement: 
let $f: A_{can} \rightarrow A_{can}$ be a homeomorphism conjugating 
the actions of $\Gamma$ on $A_{can}$ via ${\lambda}_1$, ${\lambda}_2$, respectively. 
Clearly, $f$ extends to
a homeomorphism $f: P(V) \setminus L_0 \rightarrow P(V) \setminus L_0$, still 
$\Gamma$-equivariant. 
Consider $\varphi^t$ as a group of projective transformations of $A_{can}$. 

For any nontrivial $\gamma$ in $\Gamma$, $[e_2]$ is saddle fixed point of ${\lambda}_1(\gamma)$ and ${\lambda}_2(\gamma)$, and the ${\lambda}_1(\gamma)$- or ${\lambda}_2(\gamma)$-stable leaf of $[e_2]$ is the union of two orbits of $\varphi^t$, with $[e_2]$ itself. Moreover,
when $\gamma$ is varying, these stable leaves form a dense subset among the orbits of $\varphi^t$. It follows that $f$ sends $\varphi^t$-orbits on $\varphi^t$-orbits. Hence, there is a continuous fonction $\alpha: {\mathbb R} \times A_{can} \rightarrow \mathbb R$ such that:

\[ f(\varphi^t([u])) = \varphi^{\alpha(t, [u])}f([u]) \]

For a fixed $t$, since $\varphi^t$ commutes with the $\Gamma$-actions, $[u] \mapsto \alpha(t, [u])$ is ${\lambda}_1(\Gamma)$-invariant. Since the action of ${\lambda}_1(\Gamma)$ on the annulus is minimal, it follows that $\alpha$ depends only on $t$. Since 
$f^{-1} \circ \varphi^t \circ f$ is a one-parameter subgroup of $G$, $\alpha: {\mathbb R} \rightarrow \mathbb R$ is a morphism: there is a positive constant $C$ such that $\alpha(t) = Ct$.

For some non-trivial $\gamma$ in $\Gamma$, and for $j=1,2$, denote by 
$\pm e^{a_j(\gamma)}$, $\pm e^{b_j(\gamma)}$ the eigenvalues $\neq 1$ 
of ${\lambda}_j(\gamma)$, with $b < 0 < a$. Then, the ${\lambda}_j(\gamma)$-stable 
leaf of $[e_2]$ is the fixed point set of ${\lambda}_1(\gamma)\varphi^{-b_1(\gamma)}$. 
Its image by $f$ is the fixed point set of ${\lambda}_2(\gamma)\varphi^{-Cb_1(\gamma)}$. 
It follows that $b_2(\gamma) = Cb_1(\gamma)$. Similarly, $a_2(\gamma) = Ca_1(\gamma)$. 

Consider the projections of $\lambda_j(\gamma)(\Gamma)$ in 
$\mbox{PSL}(2, \mathbb R) = \overline{H}$. They are 
Fuchsian subgroups, corresponding to metrics $g_1$, $g_2$ 
on the surface $\Sigma$ with constant curvature $-1$. Then, 
$a_2(\gamma) = Ca_1(\gamma)$ and $b_2(\gamma) = Cb_1(\gamma)$ 
imply that \emph{the $g_2$-length of a closed geodesic $c$ is $C$ 
times the $g_1$-length of the closed geodesic freely homotopic to $c$.\/} 
This is possible only if $C=1$ (see for example \cite{notes}; on page $137$, 
end of  ``Expos\'e $7$'', it is proved that the hypothesis (H) on page $136$, 
i.e., the claim $C \neq 1$, is impossible).

Therefore, if $\mbox{Tr}$ denotes the trace function on the algebra $\mbox{gl}_2$ of $2 \times 2$ matrices: 

\begin{eqnarray}
\mbox{Tr} \circ \lambda_1 = \mbox{Tr} \circ \lambda_2
\end{eqnarray}

It is well-known that this implies that $\lambda_1$ and $\lambda_2$ are conjugate in 
$\mbox{GL}_2$. Let's recall the argument: since $\lambda_j$ are irreducible 
representations, 
every element $g$ of $\mbox{gl}_2$ can be (non-uniquely) 
written as a sum $\sum_i \eta_i\lambda_j(\gamma_i)$. 

For every $g$ in $\mbox{gl}_2$, select such a decomposition 
$g = \sum_i \eta_i\lambda_1(\gamma_i)$. Then define $\phi(g) = \sum_i \eta_i\lambda_2(\gamma_i)$. The key point is that $\phi(g)$ does not depend on the selected decomposition of $g$. Indeed, if $0 = \sum_i \eta_i\lambda_1(\gamma_i)$, 
then, for every $g' = \sum_k \nu_k\lambda_2(\gamma_k)$ in $\mbox{gl}_2$, we have:

\begin{eqnarray*}
 \mbox{Tr}((\sum_i \eta_i\lambda_2(\gamma_i))g')  & = & \mbox{Tr}(\sum_{i,k} \eta_i\nu_k\lambda_2(\gamma_i)\lambda_2(\gamma_k)) \\
& = & \sum_{i,k} \eta_i\nu_k\mbox{Tr}(\lambda_2(\gamma_i\gamma_k)) \\
& = & \sum_{i,k} \eta_i\nu_k\mbox{Tr}(\lambda_1(\gamma_i\gamma_k)) \\
& = & \mbox{Tr}((\sum_i \eta_i\lambda_1(\gamma_i))(\sum_k \nu_k\lambda_1(\gamma_k))) = 0
\end{eqnarray*}

Since this holds for every $g'$, and since $g' \mapsto \mbox{Tr}(g_0g')$ can be zero only if $g_0=0$, we obtain that $0 = \sum_i \eta_i\lambda_2(\gamma_i)$. As a corollary, $\phi$ is well-defined.

This map $\phi$ is obviously an algebra automorphism of $\mbox{gl}_2$; but 
such an  automorphism is known to be an inner automorphism. Moreover, 
$\phi \circ \lambda_1 = \lambda_2$. The proposition follows.\fin

\subsection{Goldman manifolds}

Consider a $u$-linear deformation $\rho$ of a canonical Anosov flag representation $\rho_0$.
Here, we consider that the representations $\Gamma \rightarrow \mbox{GL}_2$ 
associated to $\rho$ and $\rho_0$ are the same, i.e., in the notation of the beginning of \S~\ref{flagoldman}, $\rho_0(\gamma) = e^{u(\gamma)}\rho_\lambda$.

\begin{prop}
\label{conjuguons}
The actions of $\rho(\Gamma)$ and $\rho_0(\Gamma)$ on $P(V)$ are topologically conjugate.
\end{prop}

\preu
We reproduce the short proof in \cite{barflag}: there is a (H\"older) 
continuous and homogeneous degree one map $\delta: {\mathbb R}^2 \rightarrow \mathbb R$ such that $[xe_1 + ye_3 + ze_2]$ belongs to the  invariant curve $L$ if and only if $z = \delta(x,y)$. 
The $\rho(\gamma)$-invariance of $L$ implies:

\[ \delta(\rho(\gamma)) = \delta(\gamma) + \tau_1(\gamma)x + \tau_2(\gamma)y \]

where $\tau_1(\gamma) = e^{u(\gamma)}\mu(\gamma)$, $\tau_2(\gamma) = e^{u(\gamma)}\nu(\gamma)$ (for $\mu(\gamma)$, $\nu(\gamma)$ as in \S~\ref{flagoldman}).
It follows immediately that $(x, y, z) \mapsto (x, y, z + \delta(x,y))$ induces on $P(V)$ the required topological conjugacy between $\rho_0(\Gamma)$ and $\rho(\Gamma)$.\fin

This proposition implies Theorem $5.1$ of \cite{barflag}: \emph{the first tautological foliations associated to $\rho$, $\rho_0$ are topologically conjugate.\/} It follows also that \emph{the action of $\rho(\Gamma)$ on $P(V) \setminus (L \cup [e_2])$ is minimal and uniquely ergodic.\/}

\subsection{Non hyperconvex Anosov flag manifolds.}
\label{pashyper}
Anosov flag manifolds with hyperconvex (i.e., quasi-Fuchsian) representations are fairly well understood (see remarks~\ref{fauxlorentz}, \ref{fuchspropre}). From now, we exclude this case.

\begin{lem}
\label{topodisc}
The complement of $L$ in $P(V)$ is a topological disc.
\end{lem}

\preu
If $L$ is a projective line, the lemma is obvious. Assume that $L^\ast$ 
is a projective line, i.e., that $\rho$ is a hyperbolic representation: 
consider an affine domain $U \subset P(V)$ not containing the 
$\rho(\Gamma)$-fixed point $[(L^\ast)^\perp]$. Select a coordinate 
system on $U \approx {\mathbb R}^2$ such that the intersection with $U$ 
of projective lines in $L^\ast$ are vertical lines 
$\{ \ast \} \times {\mathbb R}$. Then, according to 
corollary~\ref{lemflag1}, $L \cap U$ is the graph of a 
continuous map from $\mathbb R$ into $\mathbb R$. 
Then, the lemma becomes obvious.

We are left with the case where $L$, $L^\ast$ are not projective lines, i.e., the strongly irreducible case (cf. proposition~\ref{pointfixe?}). Let $S(V)$ be the $N$-sphere of $V$: the radial projection $\Pi_S: S(V) \rightarrow P(V)$ is the double covering, with covering automorphisms $\pm id$. The group $\mbox{SL}(V)$ acts naturally on $S(V)$: consider the sphere $S(V)$ as the space of half-lines in $V$. 

The limit curve $L$ is a simple closed curve. There are two cases:

-- either $\Pi_S^{-1}(L)$ is a connected simple closed curve $\hat{L}$,

-- or $\Pi_S^{-1}(L)$ is the union of two connected simple closed curves $\hat{L}_+$, 
$\hat{L}_-$.

In the first case, the lemma follows from Jordan's Theorem, and the $(-id)$-invariance 
of $\hat{L}$. Hence, the proof of the lemma amounts to excluding the second case, 
more precisely, to proving that in the second case, $\rho$ is hyperconvex.

The antipody $-id$ exchanges $\hat{L}_+$ and $\hat{L}_-$. For any element 
$\gamma$ of $\Gamma$, there is one and only one way to lift $\rho(\gamma)$ to an 
element $\hat{\rho}(\gamma)$ of $\mbox{SL}(V)$ preserving $\hat{L}_+$. 
It provides a representation $\hat{\rho}: \Gamma \rightarrow \mbox{SL}(V)$. 

For any $[u^\ast]$ in $L^\ast$, according to corollary~\ref{lemflag1}, 
the kernel of $u^\ast \in V^\ast$ intersects $\hat{L}_+ \cup \hat{L}_-$ in 
two points, one opposite the other. Since $\hat{L}_- = -\hat{L}_+$, the kernel 
of $u^\ast$ intersects $\hat{L}_+$ in one and only one point $u([u^\ast])$ of 
$\hat{L}_+$: hence, the sign of $u^\ast$ on $\hat{L}_+ \setminus \{ u([u^\ast]) \}$ 
is constant. Select $u^\ast \in S(V^\ast)$ so that this sign is positive.

This process defines a way to simultaneously lift $L$ and $L^\ast$ in $V$, $V^\ast$ to 
closed subsets $\hat{L}_+$, $\hat{L}_+^\ast$ such that for any $(u, v^\ast)$ 
in $\hat{L}_+ \times \hat{L}^\ast_+$ we have $\langle u \mid v^\ast \rangle \geq 0$. 
In the terminology of \cite{benoist}, this means that $\Lambda$ is a \emph{positive} 
subset of $X$. According to proposition $1.2$ of \cite{benoist}, since $\rho$ is 
strongly irreducible, $\rho(\Gamma)$ preserves a strictly convex domain of $P(V)$, 
meaning precisely that $\rho$ is hyperconvex, a contradiction.\fin

\begin{cor}
\label{comptriviale}
Every non-hyperconvex Anosov flag representation can be continuously deformed to the trivial representation.
\end{cor}

\preu
According to \cite{guichard}, a non-hyperconvex Anosov representation does not belong 
to the Hitchin component. On the other hand, it lifts to a representation 
$\hat{\rho}: \Gamma \rightarrow P^+\mbox{GL}(V)$. Indeed, keeping the notation used 
in the proof of Lemma~\ref{topodisc}, we can select a connected component $U^+$ of 
$S(V) \setminus \hat{L}$. We then define $\hat{\rho}(\gamma)$ as the unique lift 
in $P^+\mbox{GL}(V)$ of $\rho(\gamma)$ preserving $U^+$. It provides the required 
representation $\hat{\rho}$. The corollary then follows from \S~\ref{espacedef}.\fin

\begin{cor}
\label{tout}
For every $[u]$ in $P(V)$ there is an element $[v^\ast]$ in $L^\ast$ ``containing" $[u]$, i.e., such that $\langle u \mid v^\ast \rangle = 0$.
\end{cor}

\preu
If not, there is an element $u_0$ of $V \setminus \{ 0 \}$ such that $L^\ast$ 
is disjoint from $[u_0^\perp]$. The projective line $[u_0^\perp]$ then 
lifts to a great circle in $S(V^\ast)$ avoiding $\hat{L}^\ast$. Hence, $\hat{L}^\ast$ 
is contained in a hemisphere $U$. We obtain a contradiction since 
$\hat{L}^\ast = -\hat{L}^\ast$ and $U \cap (-U) = \emptyset$.\fin

Compare the following proposition with remark~\ref{fauxlorentz}:

\begin{prop}
\label{toresol}
$T$ and $T^\ast$ are topological Klein bottles, and $A$, $A^\ast$ are M\"{o}bius bands.
The domain $\Omega$ is connected, homeomorphic to the solid torus $\Omega_0 \approx {\mathbb S}^1 \times {\mathbb R}^2$.
\end{prop}

\preu
By the Jordan-Sch\"{o}nflies Theorem, there is a homeomorphism $f$ of $P(V)$ 
mapping the simple closed curve $L$ on the projective line $L_0$. 
The circle bundle $\Pi: X \rightarrow P(V)$ and pull-back bundle $f^\ast\Pi$ 
have the same Euler class; hence, $f$ lifts to some homeomorphism $F$ of $X$ into itself, 
preserving the fibers of $\Pi$, and inducing $f$ on $P(V)$. Then, $F(T)$ is the Klein 
bottle $T_0$. The first part of the proposition follows. 
The complement $W_0$ of $T_0$ in $P(V)$ is a solid torus, and $F(T^\ast) \cap  W_0$ is 
an annulus. For a better visualisation, lift to the $4$-sheeted covering 
$\widehat{X} \subset S(V) \times S(V^\ast)$: we are led to the well-known 
fact that the complement, in a compact solid torus with boundary $\widehat{W}_0$, 
of a compact annulus with boundary topologically embedded admitting as boundary 
two disjoint essential curves of $\partial\widehat{W}_0$, is the union of two solid 
tori.\fin

\begin{cor}
The flag manifold $M = \rho(\Gamma)\backslash\Omega$ is homeomorphic to a circle bundle over $\Sigma$.
\end{cor}

\preu
According to proposition~\ref{toresol}, if $\widehat{\Gamma}$ denote the fundamental group of $M$, we have an exact sequence:

\[ 0 \rightarrow {\mathbb Z} \rightarrow \widehat{\Gamma} \rightarrow \Gamma \rightarrow 0 \]
where $\mathbb Z$ is the fundamental group of the solid torus $\Omega$. 
Hence, $M$ is sufficiently large 
(its first homology group is infinite), is ${\mathbb R}P^2$-irreducible
(its universal covering is homeomorphic to ${\mathbb R}^3$), 
and is homotopically equivalent 
to a circle bundle over $\Sigma$. 
The corollary follows from \cite{topologie}.\fin

\subsection{Minimality of tautological foliations}

Recall that a $1$-dimensional foliation is \emph{minimal\/} if all the leaves are dense.

\begin{thm}
Let $\rho$ be a $(G,Y)$-Anosov representation, which is not hyperconvex. 

-- If $\rho$ is a canonical flag Anosov representation, then the tautological foliations are both minimal,

-- If $\rho$ is a a linear $u$-deformation, with $u \neq 0$, then the first tautological foliation is minimal, but the second tautological foliation is not minimal,

-- If $\rho$ is strongly irreducible, then the tautological foliations are both non-minimal.
\end{thm}

{\noindent \mbox{\bf Sketch of proof  }}
The first case follows from Proposition~\ref{conjuguons}. The first statement 
in the second case follows from Proposition~\ref{conjuguons} too. 
The second statement of this case is Theorem C of \cite{barflag}. 
The third case can be proved using the arguments of \S $5.2$ of \cite{barflag}: 
consider the set of flags $([u], [u^\ast]) \in \Omega$ such that $[u^\ast] \cap L$ 
is reduced to a point. Its closure $\mathcal M$ is invariant by 
the second tautological flow. Arguments of Lemma $5.7$ of \cite{barflag} 
proves that $\mathcal M$ is not all of $\Omega$, and those of Lemma $5.8$ 
of \cite{barflag} show that if $\mathcal M$ were empty, then the second tautological 
foliation would be expansive: by \cite{brunella}, it would be topologically
conjugate to a finite covering of some geodesic flow. This leads to a contradiction, 
as in \cite{barflag}. \fin

\section{The invariant M\"{o}bius bands}
In this section, we consider a \emph{non hyperconvex\/} Anosov flag representation $\rho$. The action of $\rho(\Gamma)$ on $\Lambda$ is conjugate to the usual projective action on ${\mathbb R}P^1$. The action on $\Omega$ is proper, topologically conjugate to the action of $\Gamma$ on $H$ by left translations, at least when $\rho$ is a deformation of a canonical flag representation inside the space of Anosov flag representations (see question $2$ in \S~\ref{conclu}).  

But the dynamics on the invariant M\"{o}bius bands are much more subtle. When $L$ is a projective line $[(u_0^\ast)^\perp]$, it follows from corollary~\ref{lemflag1} that $[u_0^\ast]$ does not belong to $L^\ast$. Hence, the flags $([u], [u_0^\ast])$ with $[u] \in L$ form a continuous $\rho(\Gamma)$-invariant curve in $A$. Conversely:

\begin{prop}
\label{measure}
Consider the Fuchsian projective action of $\Gamma \subset H$ on 
the projective line ${\mathbb R}P^1$ and the $\Gamma$-action on 
$A$ induced by $\rho$. If there is a $\Gamma$-equivariant measurable 
map $\sigma: {\mathbb R}P^1 \rightarrow A$, then the limit curve $L$ 
is a projective line. 
\end{prop}

\preu
This is essentially the content of Lemma $4.17$ of \cite{barflag}, that we reproduce 
here: let $\overline{\mathcal Q}$ be the complement of the diagonal in 
${\mathbb R}P^1 \times {\mathbb R}P^1$ (see \S~\ref{bifoliated}). 
The diagonal action of $\Gamma$ on $\overline{\mathcal Q}$ is ergodic 
for some $\Gamma$-invariant measure equivalent to the Lebesgue measure 
and preserved by the flip map $(x,y) \mapsto (y,x)$ (for example, 
the projection on the orbit space of the geodesic flow of the Liouville measure). 
We say that a subset of $\overline{\mathcal Q}$ is conull if the Lebesgue measure 
of its complement is zero. 
The crucial and classical observation is that this ergodicity property
implies that there is no measurable equivariant map from $\overline{\mathcal Q}$ into 
a topological space on which $\Gamma$ acts freely and properly discontinuously.

Decompose $\sigma: {\mathbb R}P^1 \rightarrow A$ into two functions 
$\eta: {\mathbb R}P^1 \rightarrow L$ and 
$\eta^\ast: {\mathbb R}P^1 \rightarrow P(V^\ast) \setminus L^\ast$.

Assume that the set of pairs $(x,y)$ for which
$\eta(y)$ does not belong to $\eta^\ast(x)$ is conull. Then, its
intersection with its image by the flip map 
is conull, and its intersection with all its $\Gamma$-iterates is too.
Thus, there is a conull $\Gamma$-invariant subset $\mathcal E$ of $\overline{\mathcal Q}$
of pairs $(x,y)$ for which the projective lines
$\eta^\ast(x)$ and $\eta^\ast(y)$ intersect at some point $[u(x,y)]$
different from $\eta(x)$ and $\eta(y)$.
We have then two cases: either almost every $[u(x,y)]$ belongs
to $L$, or almost all of them belong to $P(V) \setminus L$.
In the first case, 
$(x,y) \mapsto (\eta_+^{-1}(\eta(x)), \eta_+^{-1}(\eta(y)), \eta_+^{-1}([u(x,y)]))$ is a 
$\Gamma$-equivariant map
from $\mathcal E$ into the set of distinct triples of points of ${\mathbb R}P^{1}$.
Since the action of $\Gamma$ on this set of triples is free and properly
discontinuous, we obtain a contradiction with the ergodic argument
discussed above. In the second case, the map associating to a pair
$(x,y)$ the flag $([u(x,y)], \eta^\ast(x))$ is
a $\Gamma$-equivariant map from $\mathcal E$ into $\Omega$. We obtain once more a contradiction with the ergodic argument by Theorem \ref{actiopropre}.

Therefore, the measure of the set of pairs $(x,y)$
for which the line $\eta^\ast(x)$ contains $\eta(y)$ is conull. 
Then, by Fubini's Theorem, for almost every $x$ in ${\mathbb R}P^{1}$
and for almost all $y$ in ${\mathbb R}P^{1}$, $\eta(y)$
belongs to $\eta^\ast(x)$. For such a $x$, if there exists an element $\gamma$ of $\Gamma$
such that $\eta^\ast(x) \neq \rho^\ast(\gamma)\eta^\ast(x)$, then for 
almost every $y$, $\eta(y)$ is 
the unique point of intersection $[u_0]$ between $\eta^\ast(x)$ and 
$\rho^\ast(\gamma)\eta^\ast(x)$.
Hence, $[u_0]$ must be a common fixed point of all $\rho(\Gamma)$, 
but this is impossible since 
it belongs to $L$. 

It follows that $\eta^\ast(x) = \rho^\ast(\gamma)\eta^\ast(x)$ for every $\gamma$:
there is a $\rho(\Gamma)$-invariant projective line in $P(V)$. Hence, $\rho$ is hyperbolic,
and $L$ is the projective line $\eta^\ast(x)$.
\fin

Of course, the similar lemma with $A$ replaced by $A^\ast$ is true.

\begin{cor}
\label{pasregulier}
If $\rho$ is strongly irreducible, the maps $\eta_{\pm}$ are not Lipschitz.
\end{cor}

\preu
Assume that $\eta_+$ is Lipschitz. Then, it is differentiable almost 
everywhere. Its differential defines a measurable map 
$\sigma = (\eta, \eta^\ast): {\mathbb R}P^1 \rightarrow X$, 
where the first component $\eta$ is $\eta_+$, and the second 
component $\eta^\ast$ is the projective line tangent to the image 
of the differential of $\eta_+$. According to Proposition~\ref{measure},
the second component $\eta^\ast(x)$ is an element 
of $L^\ast$ for almost every $x$. According to Corollary~\ref{lemflag1}, we have 
$\eta^\ast(x) = \eta_-(x)$. More precisely, it follows from 
this corollary~\ref{lemflag1} that $L$ is locally strictly convex. 
We obtain a contradiction since $\rho$ was assumed to be non-hyperconvex.

The same proof applies to $\eta_-$.\fin

\section{Open questions}
\label{conclu}

Let $\rho:\Gamma \rightarrow G$ be a non-hyperconvex Anosov flag 
representation. Here, we have essentially answered Question $8$ of 
\cite{barflag} by the Theorem~\ref{actiopropre}. We also answered 
Question $6$: the curve $\Lambda$ is H\"{o}lder continuous. 
But some interesting questions are still open:

{\noindent \mbox{\bf Question 1}} 
\begin{em}
Is the circle bundle $\rho(\Gamma)\backslash\Omega$ homeomorphic to $\Gamma\backslash{H}$, i.e., to the double covering of the unit tangent bundle of $\Sigma = \Gamma\backslash{\mathbb H}^2$?
\end{em}

According to \cite{guichard}, the space of hyperconvex flag representations is connected: it is the entire Hitchin component. Thus, we can wonder:

{\noindent \mbox{\bf Question 2}} 
\begin{em}
Is the space of non-hyperconvex Anosov flag representations connected?
\end{em}

If the answer to question $2$ is yes, then the answer to question $1$ is yes.

{\noindent \mbox{\bf Question 3}} 
\begin{em}
Let $\rho':\Gamma \rightarrow G$ be another non-hyperconvex Anosov flag representation. Assume that the $\rho'(\Gamma)$-action on the invariant M\"{o}bius band $A'$ (respectively $(A')^\ast$) is topologically conjugate to the $\rho(\Gamma)$-action on $A$ (resp. $A^\ast$). Does that imply that $\rho'$ and $\rho$ are conjugate in $G$?
\end{em}

Observe that according to Proposition~\ref{rigid}, 
Proposition~\ref{conjuguons} and Proposition~\ref{measure}, 
the answer to this question is yes when $\rho$ is a hyperbolic radial 
representation!

{\noindent \mbox{\bf Question 4}} 
\begin{em}
Are the tautological foliations associated to non-hyperconvex Anosov flag representations
uniquely ergodic? 
\end{em}

Concerning dynamical properties of tautological foliations, recall Questions $1$ and $5$ of \cite{barflag}:

{\noindent \mbox{\bf Question 5}} 
\begin{em}
Can a tautological foliation associated to a non-hyperconvex Anosov flag representation admit a periodic orbit? Can it have non-zero entropy?
\end{em}

\section*{Appendix A: the $(G,Y)$-Anosov property for linear deformations.}
\label{AppA}

In this Appendix, we prove that canonical representations 
and their $u$-deforma\-tions are $(G,Y)$-Anosov. 
We first consider canonical representations. The proof of 
Theorem~\ref{cestanosov} we produce here is quite sophisticated, 
but it is a necessary  preparation for the most delicate case of 
$u$-deformations.

\subsection{Canonical representations}

Consider the following map $f: H \rightarrow Y$: 
$f(h) = ([\rho_0(h)e_1], [e_2], [\rho_0(h)e_3])$. 
It is a $\rho_0(\Gamma)$-equivariant map, 
defining a section $s$ of $\pi_\rho$. Our task is to prove that 
$s$ defines a $(G,Y)$-Anosov structure.

Consider $f_+=\pi_+ \circ f$: it provides a section $s^+$ of 
$\pi^X_{\rho_0}$. It can be written: 
$f_+(h) = ([\rho_0(h)e_1], [\rho_0^\ast(h)e^\ast_3])$. 

We first consider the first component $[\rho_0(h)e_1]$, 
which defines a section $s^+_P$ of the flat $P(V)$-bundle associated to 
$\rho_0$: this bundle $\pi_P: E_{\rho_0}(P) \rightarrow M$ 
is defined in the same way as the bundles $E_{\rho_0}$ and $E_{\rho_0}(X)$; 
it also admits a horizontal flow $\Phi^t_P$ above $\Phi^t$. 

Define $\Psi(h, (\alpha, \beta)) = (h, [\rho_0(h)e_1 + \alpha\rho_0(h)e_2 + \beta\rho_0(h)e_3])$. 
It is a $\Gamma$-equivariant map, when $H \times {\mathbb R}^2$ is equipped with the $\Gamma$-action $\gamma(h, (\alpha, \beta)) = (\gamma h, (\alpha, \beta))$, and $H \times P(V)$ is equipped with the action $\gamma(h, [u]) = (\gamma h, [\rho_0(\gamma)u])$.
Hence, it induces a fibered map $\bar{\Psi}: M \times {\mathbb R}^2 \rightarrow E_{\rho_0}(P)$. More precisely, it provides a trivialisation of an open neighborhood $W$ of the image of $s^+_P$.

Select any left invariant metric $m$ on $H$. Equip ${\mathbb R}^2$ with the euclidean norm $d\alpha^2+d\beta^2$. The image by $\Psi$ of the product metric is a $\Gamma$-invariant metric: it provides a metric on the neighborhood $W$. Now, we simply observe that the flow $\Phi^t_P$ is expressed in the chart $\bar{\Psi}$ by the simple expression:

\[ \Phi^t_P(p, (\alpha, \beta)) = (\Phi^t(p), (e^{t}\alpha, e^{2t}\beta)) \]

It immediately follows that the image of $s^+_P$ is a (exponentially) repellor of $\Phi^t_P$.

Similar reasoning on the second component $[\rho^\ast_0(h)e^\ast_3]$ of $f_+(h)$ shows that it provides a section of the flat $P(V^\ast)$-bundle associated to $\rho_0^\ast$ which is also a (exponentially) repellor of the corresponding horizontal flow.

Combining these two facts, we obtain that 
\emph{the image of $s^+_P$ is an (exponentially) attractor for $\Phi^t_X$.\/}

A completely similar argument establishes that the image of the section 
$s^-$ of $E_{\rho_0}(X)$ furnished by $\pi_- \circ f$ is an (exponentially) attractor for $\Phi^t_X$.
Of course, the pairs $(s^+(p), s^-(p))$ all belong to ${\mathcal Y}_{\rho_0}$. It follows that
$s \approx (s^+, s^-)$ satisfies all the criteria required in definition~\ref{def2rep}: it defines a $(G,Y)$-Anosov structure. Theorem~\ref{cestanosov} is proved.

\subsection{Linear $u$-deformations}
For a fixed inclusion ${\i}: \Gamma \subset H$ and embedding 
$\rho_0: H \rightarrow G$, any $u \in H^1(\Gamma, {\mathbb R})$ 
defines a representation $\rho_u$.
We assume here that the stable norm $\vert u \vert_s$ is less than $1/2$, and we want to prove that $\rho_u$ is $(G,Y)$-Anosov.

Since $\varphi^t$ preserves the point $[e_2]$ and acts trivially on $L = [(e_2^\ast)^\perp]$, the map $f(h) = ([\rho_0(h)e_1], [e_2], [\rho_0(h)e_3])$ is still  $\rho_u(\Gamma)$-equivariant: it provides the suitable section $s$. 

Exactly as we did for canonical representations, we study separately 
the two components $E_{\rho_u}(X)$. Each of them is then decomposed 
in bundles $E_{\rho_u}(P)$, $E_{\rho_u}(P^\ast)$.
At the end, we have to consider $4$ sections of bundles over $M$ by projective spaces, and we have to prove that the images of these sections are repellors or attractors of the associated horizontal flows.

We only discuss here the section of $E_{\rho_u}(P)$ defined by 
$h \mapsto [\rho_0(e_1)]$. The other sections can be treated in a similar way left to the reader.

We consider once more the map $\Psi(h, (\alpha, \beta)) = (h, [\rho_0(h)e_1 + \alpha\rho_0(h)e_2 + \beta\rho_0(h)e_3])$. The main difference with the canonical case is that the $\Gamma$-action on $H \times {\mathbb R}^2$ to be considered is:

\begin{eqnarray}
\gamma (h, (\alpha, \beta)) = (\gamma h, (e^{-u(\gamma)}\alpha, \beta)) 
\end{eqnarray}

The quotient of this action is then an ${\mathbb R}^2$-bundle $E$ over $M$, and we have a fibered map $\bar{\Psi}: E \rightarrow E_{\rho_u}(P)$. The image of $\bar{\Psi}$ is a neighborhood of the $\Phi^t_P$-invariant section to be studied.

The key point is to define a metric on $H \times {\mathbb R}^2$ such that:

-- the $\Gamma$-action defined by $(3)$ is isometric,

-- the null section $h \mapsto (h, 0)$ is a repellor for the horizontal flow $(h, (\alpha, \beta)) \mapsto (ha^t, (e^{t}\alpha, e^{2t}\beta))$.

Let $\Sigma$ be the Riemannian surface $\Gamma\backslash{\mathbb H}^2$. The cohomology class $u \in H^1(\Gamma, {\mathbb R})$ can be represented by a $1$-form $\omega$ on $\Sigma$ such that the integration of $\omega$ along any loop representing an element $\gamma$ of $\Gamma \approx \pi_1(\Sigma)$ is $u(\gamma)$.

The quotient $\Gamma\backslash\overline{H}$ is naturally identified with 
the unit tangent bundle of $\Sigma$: the orbit space of the left action 
of $\mbox{SO}(2)$ on $\overline{H}$ is canonically identified 
with $\Sigma$. Denote by 
$\eta: \Gamma\backslash\overline{H} \rightarrow 
\Gamma\backslash\overline{H}/\mbox{SO}(2)$ the quotient map, and consider
the $1$-form $\eta^\ast(\omega)$. The assumption $\vert u \vert_s < 1/2$ 
implies the following\footnote{Recall that $\overline{\Phi}^t$ is 
the geodesic flow with parametrisation multiplied by $2$ 
(see remark~\ref{phigeod}): this factor $2$ compensates for the 
stable norm $1/2$.\/}:
\emph{the absolute value of the integration of $\eta^\ast(\omega)$ 
along a periodic orbit of $\overline{\Phi}^t$ with period $T$ 
is less than $CT$, for $C =2\vert u \vert_s <1$.\/} 
Any orbit of $\Phi^t$ can be approximated by periodic orbits. 
Therefore, if $\bar{\theta}: [0, T] \rightarrow \overline{M}$ 
is anportion of an orbit of $\overline{\Phi}^t$ 
(i.e., $\bar{\theta}(t) = \overline{\Phi}^t(p)$ for some $p$), we have:

\[ \vert  \int_{\bar{\theta}} \eta^\ast(\omega) \vert \leq CT \]

Let $\hat{\omega}$ be the lifting of $\eta^\ast(\omega)$ to 
$M = \Gamma\backslash H$: if
$\theta(t) = \Phi^t(p)$, we have:

\begin{eqnarray}
\vert  \int_{\theta} \hat{\omega} \vert \leq CT 
\end{eqnarray}

Finally, the lifting $\tilde{\omega}$ in $H$ of $\hat{\omega}$ is exact: there is a function $\nu: H \rightarrow \mathbb R$ such that $d\nu = \tilde{\omega}$. According to $(4)$, we have:

\begin{eqnarray}
\vert \nu(ha^t) - \nu(h) \vert \leq Ct
\end{eqnarray}

Equip the fiber ${\mathbb R}^2$ over $h \in H$ with the metric $e^{2\nu(h)}d\alpha^2+d\beta^2$: it gives a metric on $H \times {\mathbb R}^2$. Since $\nu(\gamma h) - \nu(h) = \int_{[h, \gamma h]}d\nu = \int_{[h, \gamma h]} \tilde{\omega} = u(\gamma)$, the transformations defined in $(3)$ are isometries. Now, it follows from $(5)$ that the horizontal flow $(h, (\alpha, \beta)) \mapsto (ha^t, (e^{t}\alpha, e^{2t}\beta))$ expands the norm on the fibers by a factor at least $e^{(1-C)t}$. Theorem~\ref{flaganosov} follows.

\section*{Appendix B: Properness of the action in the radial case.}
\label{AppB}

In this appendix, we give a proof of \tm~\ref{actiopropre} 
in the non strongly irreducible case, i.e., according to 
remark~\ref{casreduc}, in the case of hyperbolic  
representations. We use all the notation and conventions 
introduced in \S~\ref{propre}. In this section, 
we obtained a contradiction many times by exhibiting an 
invariant projective line. This event cannot 
anymore be considered as a contradiction in the 
radial case: it slightly complicates the proof. 

On the other hand, in the hyperbolic radial case, we have of course additional 
properties:
every $g_n$ (resp. $g_n^\ast$) stabilizes $[e_2]$ (resp. $[e_2^\ast]$). 
More precisely, $e_2$, $e_2^\ast$ are eigenvectors for the 
intermediate eigenvalues $\mu_n$, $\mu_n^{-1}$. Hence, 
the $N$-isometries $k_n$, $l_n$ can be selected so that they 
both preserve $[e_2]$: $k_ne_2 = \pm e_2$, $l_ne_2 = \pm e_2$.

Let's reproduce the proof in \S~\ref{propre}, supplying the needed adaptations. Of course, Lemma~\ref{caconverge} still applies here. Let's consider the analog of Lemma~\ref{muestnul}:

\begin{lem}
\label{muestnulB}
The coefficients $\bar{\mu}$ and $\bar{\mu}'$ are both zero.
\end{lem}

\proof
Assume $\bar{\mu} \neq 0$. We prove as in \ref{muestnul} the identities:

\[ I_\ast = I^\perp, \;\;\;\;\;\; K^\perp = K_\ast \]

Moreover, $[K]$ and $[I_\ast]$ are single points.

\emph{Fact $1$: $[K]$ belongs to the invariant curve $L$.\/}
If not, the $g_n$-invariant curve $L$ is contained in $P(V) \setminus
[K]$. By  Lemma~\ref{caconverge}, it implies that $L$ is equal to the 
projective line $[I]$, which is therefore a $\rho(\Gamma)$-invariant 
projective line. In particular, $[K]$ does not belong to $[I]$. 
Let $D$ be a small open disk containing $[K]$, but with 
closure disjoint from $[I] = L$. By Lemma~\ref{caconverge}, 
the union of the iterates $g_nD$ is the complement of $[I]$ in $P(V)$. 
Hence, one of them contains the closure of $D$: it follows 
that $D$ contains the attractive fixed point of some $g_n^{-1}$. 
Taking arbitrarly small $D$, we obtain that $[K]$ is the limit of 
a sequence of attractive fixed points of elements of $\rho(\Gamma)$. 
By corollary~\ref{Lattire}, $[K]$ belongs to $L$, a contradiction.

\emph{Fact $1'$: $[I_\ast]$ belongs to $L^\ast$.\/}
Same proof as for Fact $1$, applied to the inverse sequence $(h^\ast = (g_n^\ast)^{-1})_{(n \in {\mathbb N})}$.

\emph{Fact $2$: The invariant curve $L$ is a projective line containing $[K]$.\/}
Same proof as fact $2$ in \ref{muestnul}.

\emph{Fact $2'$: The invariant curve $L^\ast$ is a projective line containing $[I_\ast]$.\/}
The proof is the same as for Fact $2$, the arguments applied to $g_n$ above to the inverse sequence $h_n^\ast =(g^\ast_n)^{-1}$.

\emph{Conclusion.\/} By facts $2$ and $2'$, the invariant curves are projective lines.
It follows that $\rho$ is actually the $u$-deformation of a canonical Anosov representation:

\[ \rho(\gamma) = \left(\begin{array}{ccc}
                  e^{u(\gamma)/3}a(\gamma) & 0 & e^{u(\gamma)/3}b(\gamma) \\
		  0 & e^{-2u(\gamma)/3} & 0 \\
		  e^{u(\gamma)/3}c(\gamma) & 0 & e^{u(\gamma)/3} d(\gamma)
\end{array}\right) \]

where $\rho_\lambda: \gamma \mapsto \left(\begin{array}{cc}
                    a(\gamma) & b(\gamma) \\
		    c(\gamma) & d(\gamma) 
		    \end{array}\right) $ is a fuchsian morphism into 
$\mbox{SL}(2,{\mathbb R})$, and $u: {\Gamma} \rightarrow {\mathbb R}$ a morphism.
According to Proposition~\ref{pointfixe?}, the stable norm of $u$ satisfies $\vert u \vert_s < 1/2$.

Lemma~\ref{muestnulB} is then one of the intermediate results of  \cite{salein}: 
indeed, let $l(\gamma) $ be as in \cite{salein} the logarithm of the spectral radius 
of the diagonal matrix appearing in the Cartan decomposition of $\gamma$ in $\mbox{SL}(2,{\mathbb R})$. Then, the logarithms of the coefficients $\lambda_n$, $\mu_n$ and $\nu_n$ are, modulo some common additive constant, the quantities 
$u(\gamma_{n}) + l(\gamma_{n})$, $u(\gamma_{n})$ and $u(\gamma_{n}) -
l(\gamma_{n})$. In the proof of Th\'eor\`eme $3.4$ in \cite{salein}, it is shown that 
$l(\gamma_{n}) - {\mid}u(\gamma_{n}){\mid}$ tends to $+\infty$: it precisely means $\bar{\mu} = 0$. \fin 

Hence, according to Lemma~\ref{muestnulB}, $[K]$ and $[K_\ast]$ are projective lines, and $[I]$, $[I_\ast]$ are single points. 

\begin{lem}
\label{IIastB}
The points $[I]$, $[I_\ast]$ belong to $L$, $L^\ast$, respectively.
\end{lem}

\preu
Assume that $[I]$ does not belong to $L$. Then, by Lemma~\ref{caconverge}, 
$L$ is contained in the projective line $[K]$. Hence, $L = [K]$. 
By considering small discs around $[I]$, we obtain, 
as in the proof of Fact $1$ in Lemma~\ref{muestnul}, that $[I]$ 
is the limit of a sequence of attractive fixed points of elements 
of $\rho({\Gamma})$. Together with Corollary~\ref{Lattire}, this 
leads to a contradiction.

A similar proof gives $[I_\ast] \in L^\ast$.\fin

\begin{lem}
\label{KKastB}
The points $[K^\perp]$ and $[K_\ast^\perp]$ belong to $L^\ast$, $L$ respectively.
\end{lem}

\preu
This is the dual version of Lemma~\ref{IIastB}: apply the reasoning above to the inverse sequence $h_n$.\fin

The proof of \tm~\ref{actiopropre} in the reducible case follows then from Lemma~\ref{caconverge} and Lemmas~\ref{IIastB}, \ref{KKastB}, exactly as the irreducible case followed in \S~\ref{propre} from Lemmas~\ref{IIast}, \ref{KKast} (and Lemma~\ref{caconverge}).

\medskip
\noindent
CNRS, UMPA, UMR 5669 \\
\'Ecole Normale Sup\'erieure de Lyon,\\
46, all\'ee d'Italie \\
69364 Lyon cedex 07, FRANCE \\

\end{document}